\documentclass[]{interact}

\usepackage{natbib}
\bibpunct[, ]{(}{)}{;}{a}{}{,}

\usepackage{url}

\usepackage[ruled,vlined,linesnumbered]{algorithm2e}

\usepackage{enumitem}
\usepackage{comment}

\usepackage{subcaption}

\usepackage{tikz,tikz-3dplot,graphicx} 
\usetikzlibrary{arrows,shadings,shadows,automata,shapes,calc,positioning,patterns,decorations.markings,decorations.pathmorphing,plotmarks,angles,backgrounds,fit,fadings,shadows.blur}

\usepackage[colorlinks = true,
            linkcolor = blue,
            urlcolor  = blue,
            citecolor = blue,
            anchorcolor = blue]{hyperref}
\newcommand{\Exp}{\mathbb{E}}
\newcommand{\Prob}{\mathbb{P}}

\newcommand{\setSce}{\mathcal{S}}
\newcommand{\Tpred}{\tau}
\newcommand{\Tpar}{T}
\newcommand{\Tset}{\mathcal{\Tpar}}

\newcommand{\Xset}{\mathcal{X}}

\newcommand{\ffnc}{f}
\newcommand{\ffnct}{\ffnc_t}

\newcommand{\xipar}{\xi}
\newcommand{\xiparrand}{{\xipar}} 
\newcommand{\xiparval}{\xipar}

\newcommand{\xisubtval}{\xiparval_t}

\newcommand{\allxi}{\text{w.p.1}}

\newcommand{\xicolor}{black} 

\newcommand{\xiT}{{\color{\xicolor}{\xiparrand_{[\Tpar]}}}}
\newcommand{\xisupt}{{\color{\xicolor}{\xiparrand_{[t]}}}}

\newcommand{\xisuptminusone}{{\color{\xicolor}{\xiparrand_{[t-1]}}}}

\newcommand{\Jpar}{J}
\newcommand{\Jset}{\mathcal{\Jpar}}

\newcommand{\ns}{N_S}
\newcommand{\np}{N_P}

\newcommand{\cred}{\color{black}}
\newcommand{\red}[1]{\textcolor{black}{#1}}

\usepackage{caption} 

\begin{document}


\title{On the Impact of Deep Learning-based Time-series Forecasts on Multistage Stochastic Programming Policies}

\author{
\name{Juyoung Wang \textsuperscript{a}, Mucahit Cevik \textsuperscript{b*}\thanks{ *Corresponding author. Email: mcevik@ryerson.ca}, and Merve Bodur \textsuperscript{a}}
\affil{\textsuperscript{a}Department of Mechanical and Industrial Engineering, University of Toronto, Toronto, ON, Canada; \textsuperscript{b}Department of Mechanical and Industrial Engineering, Ryerson University, Toronto, ON, Canada}
}

\maketitle

\begin{abstract}
Multistage stochastic programming provides a modeling framework for sequential decision-making problems that involve uncertainty. One typically overlooked aspect of this methodology is how uncertainty is incorporated into modeling. Traditionally, statistical forecasting techniques with simple forms, e.g., (first-order) autoregressive time-series models, are used to extract scenarios to be added to optimization models to represent the uncertain future. However, often times, the performance of these forecasting models are not thoroughly assessed. Motivated by the advances in probabilistic forecasting, 
we incorporate a deep learning-based 
time-series forecasting method into multistage stochastic programming framework, and compare it with the cases where a traditional forecasting method is employed to model the uncertainty. 
We assess the impact of more accurate forecasts on the quality of two commonly used look-ahead policies, a deterministic one and a two-stage one, in a rolling\red{-}horizon framework on a practical problem. 
Our results illustrate that more accurate forecasts contribute substantially to the model performance, and enable obtaining high-quality solutions even from computationally cheap heuristics. 
They also show that the probabilistic forecasting capabilities of deep learning-based methods can be especially beneficial when used as a (conditional) sampling tool for scenario-based models, and to predict the worst-case scenario for risk-averse models.
\end{abstract}






\section{Introduction}
The vast majority of mathematical programming applications assumes deterministic, static data. 
However, real world problems almost always include some \emph{uncertain parameters}. One such problem is supply chain management where product demand is highly uncertain, while electric power generation problem is another example where energy demand, water/wind inflow are possible uncertain factors at the time of generation.
It has been traditionally difficult to predict such uncertainties with high accuracy, however, with the existence of substantial historical data and advances in big data analytics, it is becoming possible to model uncertainty by fitting accurate probability distributions over uncertain parameters. 
A critical part of being able to exploit the available data is to be able to model optimization problems by taking uncertainty into account, which is the case in \emph{stochastic programming}. 
Such models provide very useful information to  decision-makers, yielding solutions that are hedged against future uncertainty. Stochastic programming has numerous applications in areas such as scheduling, production planning, supply chain management, energy, finance, manufacturing, healthcare, and natural resources \citep{wallace-ziemba:05}. 

The most widely applied stochastic programming models involve two main decision stages. In \emph{two-stage stochastic programs}, the first-stage decisions have to be made before observing random outcomes of uncertain parameters, while the second-stage decisions are made after all the uncertainty has been revealed. However, in the majority of practical problems, the planning horizon has more than two decision stages and the uncertainty is revealed gradually over time. 

\emph{Multistage stochastic programming} provides a modeling framework for sequential decision-making problems under uncertainty. 
In the classical setting, the uncertainty information is represented by a stochastic process with a given probability distribution and support, and the goal is to find a \emph{policy}, that defines decisions to be made at each decision stage, while optimizing an objective function of an expected value form. More specifically, the decisions at each stage are functions of the observed outcomes up to that stage, which should satisfy a given set of constraints almost surely with respect to the distribution of the stochastic process; and the objective is to minimize/maximize the total expected cost/profit of the feasible policies over all stages. Some risk measures other than the expected value, such as variance and conditional value at risk, can be incorporated into stochastic programs to obtain risk-averse variants, one of the most conservative one being robust optimization models that optimize the worst-case cost/profit. In this paper, we focus on risk-neutral multistage stochastic programs (MSPs), while \red{providing} one set of experiments using the robust optimization variant.

The fundamental assumption of stochastic programming is that the \emph{probability distribution} of the underlying stochastic process is given. Thus, the process of fitting an accurate (joint) probability distribution over uncertain parameters of a given decision-making problem is of utmost importance to obtain high-quality policies from a stochastic program. However, this prior process is usually glossed over in  stochastic programming applications. 
Also considering the computational difficulties associated with stochastic programs, in particular due to the famous curse of dimensionality which states the exponential growth in the problem size in terms of the number of decision stages, many crude simplifying assumptions are made for MSPs. In the literature, random variables associated with different stages are commonly assumed to be independent (referred to as stage-wise independence), and standard distribution of simple form (e.g., uniform, Gaussian, and lognormal) are fitted for random variables. In order to preserve some dependence between stages, simple (e.g., first-order) autoregressive time-series models are mostly adopted. 

Although these assumptions provide a great deal of help to relieve the computational burden for solution methods, especially for the ones that require the construction of scenario trees such as stochastic dual dynamic programming \citep{pereira1991multi}, they can be deemed inadequate to incorporate stochastic processes with high order of dependency due to the explosive state-space expansion, and in such a case, they can be even computationally intractable. 
In fact it is mostly the case that the stage-wise independence assumption is violated in practice, and the underlying probability distribution cannot be represented well with a simple form, notably as a reasonable sized scenario tree. On the other hand, \emph{deep learning-based time-series forecasting} methods can provide highly accurate models for uncertainty. Moreover, when neural network-based methods (e.g., recurrent neural networks and convolutional neural networks) are used, stage-wise dependencies can be modeled inherently. Those promising models of uncertainty can then be combined with some scenario tree-free solution methods, such as deterministic or two-stage stochastic programming approximations, and 
decision rule-based approaches, \red{which allow users to solve the optimization problem without suffering from the curse of dimensionality caused by the exponential growth of the scenario tree.} As one can expect, depending on the quality of forecasts, simpler methods as in the former class can yield high-quality solutions for MSPs, which is the focus of our computational study.


In this paper, we conduct an empirical study on the importance of the assumed probability distribution for MSPs by analyzing the impact of better forecasts. 
More specifically, we compare the quality of certain \red{look-ahead} policies under the assumption of different uncertainty models. \red{Look-ahead policies are commonly used in the context of sequential optimization problems. The term ``look-ahead" indicates that the impact of current stage decisions on the future ones are taken into account while making decisions. In the realm of sequential decision making under uncertainty, future uncertainty is also considered. There are several classes of look-ahead policies, which mainly differ in the complexity of the optimization problems and the type of forecasts for uncertain parameters.} In our work, we consider two most commonly used look-ahead policies: (i) a deterministic policy that relies on the conditional expected scenario for the future, and (ii) a two-stage policy that approximates the future uncertain stages by means of a finite set of scenarios.

We consider two commonly used primal look-ahead policies for MSPs in a rolling-horizon framework: (i) a deterministic policy that relies on the conditional expected scenario for the future, and (ii) a two-stage policy that approximates the future uncertain stages by means of a finite set of scenarios. For the scenario generation phase of these rolling-horizon schemes, we adopt a modern deep learning-based 
time-series forecasting method, namely DeepAR~\citep{Salinas2019}, and compare it with the case where a traditional 
time-series model, namely first-order autoregressive (AR) model or moving average type model, is used. In our numerical experiments on a multi-item stochastic lot-sizing problem, using two different well-known demand datasets, we find that the use of improved forecasting methods lead to significantly better policies. 
In particular, compared to AR(1), DeepAR-based solutions yield 16\% improvement on average under the risk-neutral setting based on our performance measure --- the percentage gap of the obtained objective value to the perfect information bound, \red{where the perfect information bound is defined as the expected value of the optimal values of deterministic optimization problems obtained by individual scenarios (i.e., particular realizations of the random variables in the model)}.
Moreover, in a risk-averse setting, which rely on using the prediction of the worst-case realization, DeepAR-based solutions lead to 56\% improvement on average over AR(1)-based solutions, emphasizing the enhanced probabilistic forecasting capabilities of DeepAR. 
Lastly, in the risk-neutral setting, overall, we find that two-stage policies outperform deterministic ones, and most notably DeepAR-based two-stage policies perform the best, for instance yielding 20\% improvement on the average over the AR(1)-based two-stage policies, which highlights the sampling power of DeepAR over traditional methods. All these findings reveal how helpful the use of deep learning-based forecasts can be in the realm of stochastic programming and robust optimization.

To the best of our knowledge, this is the first study incorporating modern deep learning-based time-series forecasting methods into multistage stochastic programming framework. Specifically, we make use of a recent deep-learning based method that is specifically designed to provide high-accuracy probabilistic forecasts, in arguably the most relevant context, as a sampling tool in the (multistage) stochastic optimization. Also, we illustrate its benefits in a practical setting, considering a very common application of multistage stochastic programming in supply chain management, and using real demand datasets. As such, our study contributes to bridging the gap between machine learning and operations research.


The remainder of this paper is organized as follows. In Section \ref{sec:litReview}, we review the relevant literature about MSP solution methods and time-series forecasting in the realm of stochastic programming. In Section \ref{sec:MSP}, we provide an overview of MSPs, where the necessary notation for the rest of the paper is also introduced. We present the considered policies for MSPs and the rolling-horizon framework in Section \ref{sec:Policies}, and the forecasting processes for scenario generation in Section \ref{sec:forecasting}. 
We describe our computational study setup 
in Section \ref{sec:evaluation}, and report our numerical results in Section \ref{sec:experiments}. Finally, we summarize our findings and conclude the paper in Section \ref{sec:conclusions}.

\section{Literature review}\label{sec:litReview}
Despite their expressive ability in modeling various real-life problems, MSPs have not been widely used in practice as they are notoriously difficult to solve. 
The major algorithmic challenges stem from the lack of nice properties such as convexity and continuity of the so-called value function in the existence of decisions taking only integer values, and the difficulty in computing the expectation in the objective function, especially when uncertain parameters have continuous distributions. 
As such, the majority of the literature considers the case where all decision variables are continuous, and proposes approximate solution methods, following one of the two approaches: Simplifying the underlying stochastic process to a scenario tree, and restricting the functional form of the policies. For the former case, there are a variety of available methods such as nested Benders decomposition \citep{birge1985decomposition,pereira1991multi}, aggregation and partitioning \citep{bakir2016scenario,birge1985aggregation}, and progressive hedging \citep{watson2011progressive}. 
In the latter approach, the complexity of the problem is reduced by enforcing decisions made at each stage to be a specific function of the observed outcomes up to that stage. Different functional forms, also known as decision rules, are proposed in the literature such as linear \citep{shapiro2005complexity}, piecewise linear  \citep{chen2008linear}, polynomial \citep{bampou2011scenario} and two-stage \citep{bodur2018two}. 
On the other hand, the literature is rather limited for MSPs involving integer variables. 
Regarding the scenario tree-based approaches, there exist some generic decomposition methods \citep{caroe1999dual,lulli2004branch}, as well as some algorithms designed for MSPs with binary and continuous variables only \citep{alonso2003bfc,zou2019stochastic}. For decision rules, the standard approach for the pure continuous case is extended to the binary case via piecewise linear binary functions as the policy form \citep{bertsimas2015design}, and more recently Lagrangian dual decision rules are proposed for multistage stochastic mixed integer programs \citep{daryalal2020lagrangian}.

Time-series forecasting methods are designed to take spatio-temporal relations into account. Earlier studies on time-series forecasting focus on linear prediction models such as autoregressive (AR), moving average (MA) and auto-regressive integrated moving average (ARIMA) models where a linear function of past observations is used to predict the future values \citep{Box2015}. Recent advances in artificial neural networks and deep learning allow practitioners to use deep learning in time-series forecasting, which provides the ability to process a large amount of data (e.g., a larger portion of the temporal data). Long short-term memory neural networks and gated recurrent units are among the most commonly used deep learning approaches for time-series forecasting with immediate applications in various fields such as power systems and marketing \citep{agrawal2018long, cheng2017powerlstm, kuan2017short}. Recent deep learning-based time-series prediction studies tend to involve more advanced neural network architectures such as DeepAR \citep{Salinas2019} and Deep State Space \citep{rangapuram2018deep} models. Unlike traditional time series forecasting models, such \red{deep learning-based} models leverage information from multiple time series to yield improved forecasts for any target time series. We refer readers to recent review papers \citep{Fawaz2019, Gamboa2017} on deep learning methods in time-series prediction for a more detailed list of relevant studies.

In order to overcome some computational difficulties associated with MSPs, that are exacerbated when there are many uncertain parameters in the problem, uncertainty models of a simple form have been mostly used in the literature. Under the common independence assumption, a model is built for each random variable, mostly in the form of a uniform, normal or lognormal distribution. In a more complex setting where the stage-wise independence assumption has been relaxed, previous studies mostly rely on low-order autoregressive models, namely AR(1) and AR(2), to capture some dependencies in uncertain parameters, as in
\citep{daryalal2020lagrangian,shapiro2014lectures}.

However, the accuracy of the forecasts might have a significant impact on the MSP model results, for both the actual quality of obtained solutions and accuracy of used performance measures, as well as on the number of scenarios required to obtain stable objective function values, thus on the complexity of the solution approach. 

As mentioned above, accurate inclusion of uncertainty to optimization models is a task of great importance. Some recent works address this issue and discuss the value of utilizing a precise representation of uncertainty by incorporating time-series modeling techniques relatively complex compared to the traditionally used AR(1) model. For instance, a robust optimization model for hydro-electric reservoir management problems is proposed in \citep{gauvin2018stochastic}, where the uncertainty set is built with the help of ARMA and GARCH models. It is also proven that ARMA model can be used to build an ellipsoidal type of uncertainty set to model the uncertainty of water inflow to the hydro-electric reservoir. Due to the stage-wise dependence of used time-series models, stochastic dual dynamic programming could not be used, instead affine decision rules are used to derive approximate solutions. The numerical studies in \citep{gauvin2018stochastic} indeed provide evidence of acquirable profit in employing sophisticated uncertainty modeling approaches, instead of simpler ones. 
In \citep{gauvin2018successive}, a nonlinear time-series uncertainty set modeling approach is proposed for the same problem class of \citep{gauvin2018stochastic}. Unlike the model proposed in \citep{gauvin2018stochastic}, the introduced nonlinearity of time-series makes the model non-convex. To cope with such an issue, a specialized successive linear programming algorithm is proposed. For a comprehensive review of how uncertainties are modeled in the context of energy generation planning, we refer readers to \citep{lorca2020challenges}. To the best of our knowledge, modern deep learning-based time-series forecasting techniques have not been incorporated into stochastic programming.

Quantitative studies on the benefit of precisely modeling uncertainty are performed also in contexts different than energy generation planning, such as call center demand prediction. In \citep{ye2019call}, a time-series model is developed, and it is shown that under some dependence condition of call center demand streams, performing a joint forecast is more beneficial than building a single prediction model for each stream separately. More specifically, three types of inter-stream dependence are considered and modeled simultaneously, as well as certain data conditions where the multi-stream approach performs better than the single stream model are provided. This finding also highlights the importance of taking dependencies into account while modeling uncertainty. In our approach, \red{we use }
DeepAR 
which 
considers such dependencies by design while building a model for uncertainty.   

\section{Multistage stochastic programs}\label{sec:MSP}
In this section, we introduce a generic multistage stochastic programming formulation, mostly following the presentation of \citep[Chapter 3]{shapiro2014lectures}. A generic $\Tpar$-stage stochastic programming problem, driven by an exogenous random data (i.e., stochastic) process, can be modeled as follows:
\begin{subequations}
\label{eqs:msp}
\begin{align}
	\min  \ \ & \Exp_{\xiT}\left[\sum_{t\in \Tset} \ffnct\left(x_t(\xisupt),\xisubtval\right)\right]  \label{eq:mspobj} \\
	\text{s.t.} \ \ & x_t(\xisupt)\in \Xset_t\left(x_{t-1}(\xisuptminusone),\xisubtval\right),  t\in \Tset, \allxi \label{eq:mspconsts} 
\end{align}
\end{subequations}
\noindent where
\begin{itemize}
\setlength{\itemsep}{0.1cm}
    \item $\Tset = \{ 1,2,\hdots, \Tpar \}$ is the set of decision \emph{stages}.
    \item $\xisubtval$ is the vector of random variables at stage $t \in \Tset$.
    \item $\{ \xisubtval \}_{t \in \Tset}$ represents the underlying stochastic process.
    \item $\xisupt = (\xipar_1,\xipar_2,\hdots,\xisubtval)$ is the history of the stochastic process up to stage $t \in \Tset$. 
    \item $\{x_t(\xisupt)\}_{t \in \Tset}$ describe a \emph{policy}, a solution to the MSP. Due to the underlying randomness, they are functions of the stochastic process. However, they are \emph{nonanticipative}, i.e., only depend on the history of the process up to the stage they are to be made.
    \item It is assumed that $\xipar_1 = 1$, i.e., the first stage is \emph{deterministic}.
    \item $x_0(\xipar_0)$ represents the initial state of the system.
    \item $\ffnct(\cdot,\cdot)$ is the cost function associated with stage $t \in \Tset$; it maps all the decisions made up to stage $t$, most notably the stage-$t$ decisions, and the lastly observed random outcomes (i.e., the realization of $\xisubtval$) into a real value.
    \item $\Xset_t(\cdot,\cdot)$ describes the feasible set for decisions to be made at stage $t \in \Tset$; it maps all the decisions made up to the previous stage, $t-1$, and the lastly observed outcomes, of the random vector $\xisubtval$, into a set of real-valued vectors.
\end{itemize}
This modeling framework is illustrated in Figure \ref{fig:msp}.

\begin{figure}[h]
\centering
\scalebox{0.9}{
\begin{tikzpicture}[arrow/.style={thick,->,shorten >=1pt,shorten <=1pt,>=stealth},]
	\tikzstyle{block} = [rectangle, draw, text centered, minimum height=1.8em, minimum width=0.3em]
	
	\node[block,midway,above,yshift = 0.5em] (A) {
    		 \renewcommand\arraystretch{0.5} $\begin{array}{c} \text{Stage $t-1$} \\[0.1em] \text{decisions} \end{array}$};
    		 
		 \node[right of=A,xshift = 3.3cm,block,midway,above,yshift = 0.5em] (B) {
    		 \renewcommand\arraystretch{0.5} $\begin{array}{c} \text{Stage $t$} \\[0.1em] \text{decisions} \end{array}$};
		 
	
	\draw[arrow] (A.east) to node[sloped,swap,anchor=center, above][yshift=0.5pt] {\footnotesize\parbox{1.7cm}{\centering observe\\ realization}} (B.west); 
	
	\draw[arrow] (A.east) to node[sloped,swap,anchor=center, below][yshift=-0.5pt] {\footnotesize\parbox{1.7cm}{\centering of random \\ vector $\xisubtval$}} (B.west);

       \node[left of=A,xshift = -0.7cm,midway,above,yshift = 1.2em] (D1) {$\cdots$};
       \node[right of=B,xshift = 4.8cm,midway,above,yshift = 1.2em] (D2) {$\cdots$};
       \node[left of=D1,xshift = -1.8cm,block,midway,above,yshift = 1em] (T0) {$t = 1$};
       \node[left of=A,xshift = -4.15cm,block,dashed,midway,above,yshift = 0.5em] (state0) {\renewcommand\arraystretch{0.5} $\begin{array}{c} \text{Initial} \\[0.1em] \text{State} \end{array}$};
         
         \draw[arrow,dotted] (state0) to node[sloped,swap,anchor=center, above][yshift = 0cm] {\footnotesize\text{$\xipar_1 = 1 $}} (T0);
         
       \node[right of=D2,xshift = 5.85cm,block, midway,above,yshift = 1em] (T) {$t = \Tpar$};

        \node[below of=A, midway,above,yshift = 1em,xshift = 0.15cm] (s) {$ x_{t-1}(\xisuptminusone)$};
        
         \node[below of=B,xshift = 4.3cm, midway,above,yshift = 1em] (x) {$x_t(\xisupt)$};
        

         \node[left of=x, xshift = -6.2cm] (x1) {{ $x_1(\xipar_1)$}};
         
         \node[right of=x, xshift = 1.6cm] (xT) {{ $x_T(\xiT)$}};
         
         \node[right of=x1, xshift = -3.3cm] (x0) {{ $x_0(\xipar_0)$}};
         
         \node[below of=x1, yshift = -0.4cm] (f1) {{ $\ffnc_1(x_1(\xipar_1),\xipar_1)$}};
         \draw[arrow,dotted] (x1) -- (f1);
         
         \node[below of=s, yshift = -0.2cm] (fsdummy) {};
         \node[below of=s, yshift = -0.4cm, xshift=0.2cm] (fs) {{ $\ffnc_{t-1}(x_{t-1}(\xisuptminusone),\xipar_{t-1})$}};
         \draw[arrow,dotted] (s) -- (fsdummy);
         
         \node[below of=x, yshift = -0.4cm] (fx) {{ $\ffnct(x_t(\xisupt),\xipar_{t})$}};
         \draw[arrow,dotted] (x) -- (fx);
         
         \node[below of=xT, yshift = -0.4cm] (fT) {{ $\ffnc_\Tpar(x_\Tpar(\xiT),\xipar_\Tpar)$}};
         \draw[arrow,dotted] (xT) -- (fT);
         
         \node[right of=xT, xshift = 0.68cm] (P) {$\;\;\rightarrow$ \footnotesize\text{Policy}};
         \node[right of=fT, xshift = 0.78cm] (C) {$\;\;\rightarrow$ \footnotesize\text{Cost}};
\end{tikzpicture}
}
\vskip -0.5cm
\caption{Multistage stochastic programming framework.}
\label{fig:msp}
\end{figure}

The objective function \eqref{eq:mspobj} minimizes the expected total cost of the selected policy over all decision stages. Note that the expectation is taken over the (joint) distribution of the full stochastic process. Constraints \eqref{eq:mspconsts} ensure that a feasible policy is selected, i.e., the functional decision variables at every stage satisfy the associated constraints with probability one (denoted by \allxi). We note that some decision variables might be integer, which is also enforced by \eqref{eq:mspconsts}.


As in general, we assume that the model \eqref{eqs:msp} 
\begin{itemize}
    \item is feasible and has an optimal solution (policy), and
    \item has relatively complete recourse, i.e., for any stage $t \in \Tset$, given any observation of the history up to the current stage, $\xisupt$, and any feasible set of decisions up to the previous stage, $\{  x_{t'-1}( {\color{\xicolor}{\xiparrand_{[t'-1]}}} ) \}_{t' \in \{ 1, \hdots,t-1 \} }$, there exists a feasible set of decisions for the current stage, $x_t(\xisupt)$.
\end{itemize}

Lastly, we note that although the solution to an MSP is a policy (which describes the decisions to be made at each stage as a function of the observed history up to that stage), the only \emph{implementable decisions} are the first-stage ones, i.e., $x_1(\xipar_1)$, as only the first stage is deterministic (since $\xipar_1 = 1$ is a constant, not a random variable); all the other decision variables are random variables themselves.

\section{Policies and rolling-horizon framework}\label{sec:Policies}

The model \eqref{eqs:msp} constitutes an infinite-dimensional optimization problem as both the decision variables and constraints are of infinite size, \red{as both of them are function of random variables with possibly infinite support.} Therefore, it is both theoretically and computationally very challenging to derive solutions for it, even in the case where all the objective and constraint functions are linear. The existence of integer variables adds another level of difficulty. Therefore, some heuristic procedures are preferred to construct feasible policies, or rather implementable decisions (i.e., the first-stage solution). 
More specifically, candidate first-stage solutions are usually obtained by solving a restriction or a relaxation of the original multistage model. For instance, 
\begin{itemize}
    \item a one-stage relaxation can be obtained by ignoring all the future stages, i.e., by removing all the variables $\{x_t(\xisupt)\}_{t \in \Tset \setminus \{ 1 \} }$ and the constraints associated with them from the model, and solving the remaining single-stage deterministic model; 
    \item a one-stage restriction of the multistage model can be obtained by making all the decision variables deterministic, i.e., for all $t \in \Tset \setminus \{ 1 \}$, by replacing $x_t(\xisupt)$ with $x_t(\xipar_1)$, taking away their flexibility to depend on the observed history; or by applying static linear decision rules \citep{shapiro2005complexity};
    \item a two-stage relaxation can be obtained by relaxing the nonanticipativity constraints on all the decision variables except the first-stage ones, i.e., for all $t \in \Tset \setminus \{ 1 \}$, by replacing $x_t(\xisupt)$ with $x_t(\xiT)$  to make them functions of the full history (up to $\Tpar$) despite the fact that it will not be available at stage $t$;  and
    \item a two-stage restriction can be obtained by applying two-stage linear decision rules \citep{bodur2018two}.
\end{itemize} 
Then, in practice, a full implementable solution for a $T$-stage problem can be obtained by iteratively applying the selected restriction or relaxation technique at every stage, i.e., by applying the chosen first-stage solution generation method in a \emph{rolling-horizon framework}, which is illustrated in Algorithm \ref{alg:rollinghorizon}. 
\red{Rolling-horizon evaluation framework is commonly used 
for discrete-time sequential decision-making problems under uncertainty \citep{silvente2018rolling, bischi2019rolling, daryalal2020lagrangian}.}


\begin{algorithm}
	\caption{Rolling-horizon framework}
	\label{alg:rollinghorizon}
	\KwIn{$\ns$, the number of stages, and $\np$,  the number of periods ($\np \geq \ns$), $\hat{x}_0$, vector representing the initial state of the system}  
	\KwOut{An implemented policy $\{ \hat{x}_t \}_{t \in \Tset}$}  
	\For{$t=1,2,\hdots,T$}{
	    $\hat{\xipar}_t \leftarrow$ observed $\xisubtval$ realization \\
	    $\mathcal{P} = \{ t,t+1,\hdots, \min \{ \Tpar, t + \np - 1\} \} $ \\
	    \If{$t > \Tpar - \ns + 1$}{ 
	    $\ns = T - t + 1$
	    }
	    $\{ \mathcal{P}_n \}_{n = 1,2,\hdots,\ns} \leftarrow $ \textbf{PartitionPeriods}$(\mathcal{P},\ns)$ \\
	    $\hat{x}_t \leftarrow$ \textbf{Optimize}$(\hat{x}_{t-1},\hat{\xipar}_{[t]},\{ \mathcal{P}_n \}_{n = 1,2,\hdots,\ns})$
	}
\end{algorithm}
\DecMargin{1em}

\begin{figure*}[ht!]
\centering
\scalebox{0.9}{
\begin{tikzpicture}[arrow/.style={thick,->,shorten >=1pt,shorten <=1pt,>=stealth},]
	\tikzstyle{block} = [rectangle, draw, text centered, minimum height=1.8em, minimum width=0.5em]
	
	\node[block,midway,above,yshift = 1.2em,xshift = 3cm] (A) {
    		 \renewcommand\arraystretch{0.5} $\begin{array}{c} x_4(\xipar_4,\xipar_5) \\[0.25em] x_5(\xipar_4,\xipar_5) \end{array}$};
    		 
		 \node[right of=A,xshift = 6cm,block,midway,above,yshift = 1.2em] (B) {
    		 \renewcommand\arraystretch{0.5} $\begin{array}{c} x_6(\xipar_6,\xipar_7) \\[0.25em] x_7(\xipar_6,\xipar_7) \end{array}$};
		 
	\draw[arrow] (A.east) to node[sloped,swap,anchor=center, above][yshift=0.5pt] {$\xipar_6, \xipar_7$} (B.west);
	
       \node[left of=A,xshift = 1cm,block,midway,above,yshift = 1.65em] (C) {$x_3$};
       
         
         \node[left of=A,xshift = -1cm,block,dashed,midway,above,yshift = 1.65em] (state0) {$\hat{x}_2$};
         
         \draw[arrow,dotted] (state0) to node[sloped,swap,anchor=center, above][yshift=0.5pt] {$\hat{\xipar}_{[3]}$} (C);
         
        \draw[arrow] (C.east) to node[sloped,swap,anchor=center, above][yshift=0.5pt] {$\xipar_4, \xipar_5$} (A.west);
         
    \draw[x=2cm,xshift=-6cm] (1,0) grid (8,0);
    \foreach \x in {1,2} \draw (2*\x-6,0.1) -- (2*\x-6,-0.1) node[below] {\x};
    \foreach \x in {4,5,...,8} \draw (2*\x-6,0.1) -- (2*\x-6,-0.1) node[below] {\x};
    
    \foreach \x in {3} \draw (2*\x-6,0.1) -- (2*\x-6,-0.1) node[below] {$t=\x$};
         
\end{tikzpicture}
}
\caption{Rolling horizon framework.}
\label{fig:rolling}
\end{figure*}
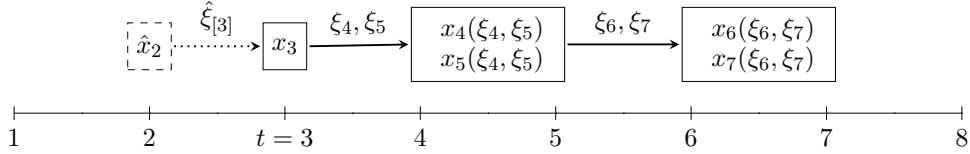

\red{In Algorithm \ref{alg:rollinghorizon}, a} policy can be constructed gradually by solving a sequence of optimization problems considering $\np$ consecutive stages ahead, each being referred to as a \emph{period}, but grouping them into $\ns$ decision \emph{stages} (which is usually small such as one or two). In the algorithm,
\begin{itemize}[leftmargin=0.5cm]
    \item PartitionPeriods($\cdot$) function takes the list of considered periods and partitions them into $\ns$ sets such that (i) each partition set consists of consecutive periods, and (ii) the current stage $t$ belongs to the first partition set $\mathcal{P}_1$; while
    \item Optimize($\cdot$) function builds and solves an $\ns$-stage stochastic program where stage $n$ consists of periods $\mathcal{P}_n$, where the initial state is fixed to $\hat{x}_{t-1}$, and returns the optimal first-stage solution. 
\end{itemize}
For instance, Figure \ref{fig:rolling} illustrates the $\ns = 3$-stage stochastic model to be solved at $t = 3$ considering the $\np = 5$ periods ahead, $\mathcal{P} = \{ 3,4,5,6,7 \}$, partitioned as $\mathcal{P}_1 = \{ 3 \}, \mathcal{P}_2 = \{ 4,5 \},\mathcal{P}_3 = \{ 6,7 \}$ for a $\Tpar = 8$-stage problem in Algorithm \ref{alg:rollinghorizon}.

In our experiments, we consider all the future stages ($\np = \Tpar$) at any stage, and use a deterministic policy ($\ns = 1$, $\mathcal{P}_1 = \mathcal{P}$) as well as a two-stage stochastic policy ($\ns = 2, \mathcal{P}_1 = \{ 1 \}, \mathcal{P}_2 = \mathcal{P} \setminus \{ 1 \}$). 

Next, we provide the details of the Optimize function for the considered policies, namely the associated mathematical models and how to solve them in conjunction with forecasting/sampling methods. We let $\Prob$ denote the joint probability distribution of all random variables involved in the stochastic process $\{ \xisubtval \}_{t \in \Tset}$. 


\subsection{Deterministic policy}
\label{subsec:det}
For deterministic policies, at stage $t \in \Tset$ of the rolling-horizon framework, given the previous stage optimal first-stage solution $\hat{x}_{t-1}$ as the initial state, and the history of observations $\hat{\xipar}_{[t]}$, we solve the following  $(\Tpar-t+1)$-period deterministic problem:
\begin{align*}
	\min  \ \ & \ffnct\left(x_{t}, \hat{\xipar}_{t} \right) + \sum_{t' = t+1}^\Tpar \ffnc_{t'}\left(x_{t'}, \xipar^{\text{expected}}_{t'} \right)  \\
	\text{s.t.} \ \ & x_t \in \Xset_{t}\left(\hat{x}_{t-1},\hat{\xipar}_t\right) \\
	& x_{t'} \in \Xset_{t'} \left(x_{t'-1},\xipar^{\text{expected}}_{t'}\right), \ t' = t+1, \hdots, \Tpar  
\end{align*}
where for each period $t' = t+1, \hdots, \Tpar$, the expected realization of the associated random vector is obtained by computing its expectation conditioned on the observed history:
\begin{equation}
\label{eq:expscen}  
\xipar^{\text{expected}}_{t'} := \Exp_{\Prob}[\xipar_{t'} | \hat{\xipar}_{[t]}]
\end{equation}
In other words, in this policy, we replace the future by the single conditional expected scenario. 

\begin{figure}[!ht]
    \centering
    \subfloat[Deterministic policy. \label{fig:deterministic}]{
    \begin{tikzpicture}[arrow/.style={thick,->,shorten >=1pt,shorten <=1pt,>=stealth},]
	\tikzstyle{block} = [rectangle, draw, text centered, minimum height=1.8em, minimum width=0.5em]
         
    \draw[x=2cm,xshift=-6cm] (1,0) grid (4,0);

    \foreach \x in {1,2,...,3} \draw (2*\x-6,0.1) -- (2*\x-6,-0.1) node[below] (T\x) {$t=\x$};
    \draw (2*4-6,0.1) -- (2*4-6,-0.1) node[below] (T4) {$\Tpar=4$};
    
    \node[above of=T4,yshift = 1.7cm,xshift = 0.1cm] (A) {};
    \node[above of=T4,yshift = 1cm,xshift = 0.1cm] (B) {};
    \node[above of=T4,yshift = 0.2cm,xshift = 0.1cm] (C) {};
    
    \draw[dotted,thick] (-4,0) -- (A);
    \draw[dotted,thick] (-2,0) -- (B);
    \draw[dotted,thick] (0,0) -- (C);
    
    \node[above of=T2,yshift = 0.4cm,xshift = -0.08cm,rotate=23] (T1E2) {\footnotesize $\Exp_{\Prob}[\xipar_{2}]$};
    \node[above of=T3,yshift = 1.15cm,xshift = -0.08cm,rotate=23] (T1E3) {\footnotesize $\Exp_{\Prob}[\xipar_{3}]$};
    \node[above of=T4,yshift = 1.9cm,xshift = -0.08cm,rotate=23] (T1E4) {\footnotesize $\Exp_{\Prob}[\xipar_{4}]$};
    
    \node[above of=T3,yshift = 0.45cm,xshift = 0cm,rotate=23] (T2E3) {\footnotesize $\Exp_{\Prob}[\xipar_{3} | \hat{\xipar}_{[2]}]$};
    \node[above of=T4,yshift = 1.25cm,xshift = 0cm,rotate=23] (T2E4) {\footnotesize $\Exp_{\Prob}[\xipar_{4} | \hat{\xipar}_{[2]}]$};
    
    \node[above of=T4,yshift = 0.45cm,xshift = 0cm,rotate=23] (T3E4) {\footnotesize $\Exp_{\Prob}[\xipar_{4} | \hat{\xipar}_{[3]}]$};
\end{tikzpicture}
    }
    \subfloat[Two-stage policy. \label{fig:twostage}]{
    \begin{tikzpicture}[arrow/.style={thick,->,shorten >=1pt,shorten <=1pt,>=stealth},]
	\tikzstyle{block} = [rectangle, draw, text centered, minimum height=1.8em, minimum width=0.5em]
         
    \draw[x=2cm,xshift=-6cm] (1,0) grid (4,0);

    \foreach \x in {1,2,...,3} \draw (2*\x-6,0.1) -- (2*\x-6,-0.1) node[below] (T\x) {$t=\x$};
    \draw (2*4-6,0.1) -- (2*4-6,-0.1) node[below] (T4) {$\Tpar=4$};
    
    \node[above of=T4,yshift = 1.6cm,xshift = 0.1cm] (A1) {};
    \node[above of=T4,yshift = 1.4cm,xshift = 0.1cm] (A2) {};
    \node[above of=T4,yshift = 1.2cm,xshift = 0.1cm] (A3) {};
    
    \node[above of=T4,yshift = 0.8cm,xshift = 0.1cm] (B1) {};
    \node[above of=T4,yshift = 0.6cm,xshift = 0.1cm] (B2) {};
    \node[above of=T4,yshift = 0.4cm,xshift = 0.1cm] (B3) {};
    
    \node[above of=T4,yshift = -0.2cm,xshift = 0.1cm] (C1) {};
    \node[above of=T4,yshift = -0.4cm,xshift = 0.1cm] (C2) {};
    \node[above of=T4,yshift = -0cm,xshift = 0.1cm] (C3) {};
    
    \draw[dotted,thick] (-4,0) -- (A1);
    \draw[dotted,thick] (-4,0) -- (A2);
    \draw[dotted,thick] (-4,0) -- (A3);
    \draw[dotted,thick] (-2,0) -- (B1);
    \draw[dotted,thick] (-2,0) -- (B2);
    \draw[dotted,thick] (-2,0) -- (B3);
    \draw[dotted,thick] (0,0) -- (C1);
    \draw[dotted,thick] (0,0) -- (C2);
    \draw[dotted,thick] (0,0) -- (C3);
    
    
    
\end{tikzpicture}
    }
    \caption{Look-ahead policies.}
    \label{fig:policies}
\end{figure}
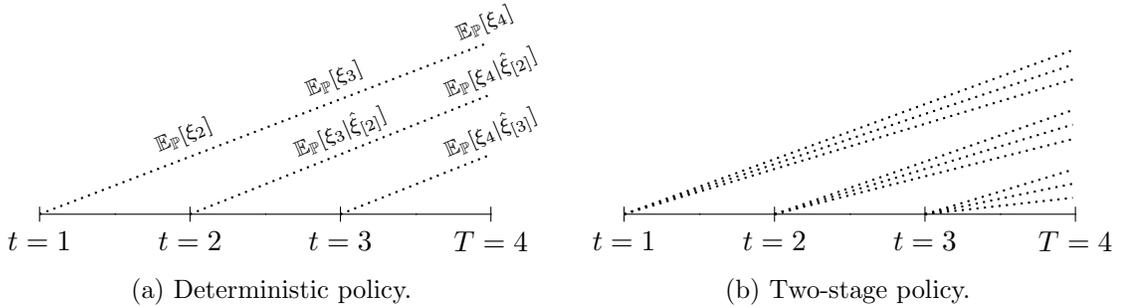

\red{An illustration for this policy is provided in Figure \ref{fig:deterministic}. As shown in the figure, the deterministic policy does not rely on multiple sample paths to describe the uncertain nature of the problem. Instead, at each decision-making stage, the uncertain future is represented by a single sample path obtained by using the conditional mean values of the random variables. As a result, such a policy is computationally cheap. However, as the uncertain factors are replaced by a single future scenario, the variance of the stochastic process is not taken into account. 
To address such an issue, a two-stage policy that approximates future uncertainty by a set of scenarios, as illustrated in Figure \ref{fig:twostage}, can be considered as an alternative, which is explained in more detail next.} 

\subsection{Two-stage policy}
\label{subsec:stoch}
In the case of two-stage policies, the future is represented by a finite set of scenarios rather than a single scenario. 
Fix a stage $t \in \Tset$ of the rolling-horizon framework. 
\begin{equation}
\label{eq:scens} 
\{ (\xipar^s_{t+1},\hdots,\xipar^s_{\Tpar} ) \}_{s \in \setSce}
\end{equation}
be an independent and identically distributed sample of future scenarios generated using the distribution $\Prob$ conditioned on the observed history $\hat{\xipar}_{[t]}$\red{, where $\mathcal{S}$ is the set of scenarios and $s$ is the scenario index.} Then, we solve the following two-stage stochastic program: 
\begin{align*}
	\min  \ \ & \ffnct\left(x_{t}, \hat{\xipar}_{t} \right) + \frac{1}{|\setSce|}\sum_{s \in \setSce} \sum_{t' = t+1}^\Tpar \ffnc_{t'}\left(x^s_{t'}, \xipar^{\text{s}}_{t'} \right)   \\
	\text{s.t.} \ \ & x_t \in \Xset_{t}\left(\hat{x}_{t-1},\hat{\xipar}_t\right) \\
	& x^s_{t'} \in \Xset_{t'} \left(x^s_{t'-1},\xipar^{\text{s}}_{t'}\right), \ t' = t+1, \hdots, \Tpar, s \in \setSce  
\end{align*}
where the objective function minimizes the average total cost over the future scenarios. 
\red{
Unlike the only implementable non-anticipative decisions $x_t$, we generate the scenario copies of future decisions, $x^s_{t'}$.
}
As this model does not scale well computationally with the sample size, we usually represent the future using a small set of scenarios.

\section{Forecasting for scenario generation}\label{sec:forecasting}
In this section, 
we describe our proposed approach, followed by the review of forecasting methods employed in our analysis. 

\subsection{Proposed approach}
As the future is represented by only a single scenario, the deterministic approach usually provides a very crude approximation of the underlying MSP, thus yields poor-quality policies. Although the two-stage policy overcomes this issue in theory, it may not be a viable approach in practice either, as it usually requires a large number of scenarios to get good-quality solutions, thus possibly a prohibitive computational effort. Furthermore, if the assumed probability distribution $\Prob$ is far from the true distribution, one can not expect the two-stage policy perform well even if a large number of scenarios have been considered. Therefore, predicting $\Prob$ accurately is of paramount importance to obtain high-quality policies. Inspired by the success of recent deep learning-based time-series forecasting methods \citep{Gamboa2017, rangapuram2018deep, Salinas2019}, we propose to estimate $\Prob$ by such a modern probabilistic forecasting method, to be given as an input to an MSP. In Section \ref{sec:experiments}, we showcase the impact of improved forecasting methods on MSP policy generation via our empirical study. 

\subsection{Forecasting methodologies}
We compare the quality of the deterministic and two-stage policies using three forecasting methods: AR(1), (log) moving average and DeepAR. We use (log) moving average as a middle-level method in the sense that it is sophisticated enough compared to AR(1), while it is still primitive compared to state-of-the-art deep learning-based time-series forecasting methodologies, such as DeepAR \citep{Salinas2019} and Deep State Space \citep{rangapuram2018deep} models, among others. Also, MA-based model has been shown to perform well for one of the datasets that we use in our experiments.

\subsubsection{AR(1)}
AR(1) is one of the most traditionally used forecasting methods in the realm of stochastic programming, due to its theoretical suitability to be applied in many existing MSP solution methodologies. 
In this modeling scheme, stochastic process $\{\xipar_t\}_{t \in \Tset}$ is assumed to have the form
\begin{align*}
\xipar_t = \phi \xipar_{t-1} + (1-\phi) \gamma + \varepsilon_t, \ \text{for all} \ t =2,\hdots,\Tpar
\end{align*}
where $\phi$ is the AR(1) coefficient, $\gamma$ is the trend constant vector, and $\varepsilon_t$ is assumed to be a white noise vector. This model
is widely used in the stochastic dual dynamic programming \citep{pereira1991multi} literature, due to its suitability for state-space expansion. For further details regarding the use of AR(1) model in the MSP context, we refer readers to \citep{dowsonsddp, shapiro2011analysis}. 


\subsubsection{Autoregressive moving average (ARMA)}
ARMA, parametrized by two natural numbers $p$ and $q$, is commonly used \citep{gauvin2018stochastic, shapiro:ejor13} forecasting methods that linearly models serial dependency of discrete time dependent sequence of random variables \citep{Box2015}, where the randomness is represented by a white noise term. A stochastic process $\{ \xipar_t \}_{t \in \Tset}$ following ARMA$(p,q)$ can be written in the following manner:\\
\begin{align*}
\xipar_t = c + \varepsilon_t + \sum_{i=1}^p \varphi_i \xipar_{t-i} + \sum _{i=1}^q \red{\psi}_i \varepsilon_{t-i}, \ \text{for all} \ t =2,\hdots,\Tpar
\end{align*}
where $\{\varepsilon_t\}_t$ is a sequence white noise, while $\{\varphi_i\}_i$ and $\{\red{\psi}_i\}_i$ are coefficients used to fit the model. The $\xipar_{t'}$ and $\varepsilon_{t'}$ values for $t' \leq 0$ are either known non-zero values or assumed to be zero. Note that by restricting $(p,q) = (1,0)$, we recover the AR(1) model. 

In our computational study, for one of the datasets, we use a logarithmic moving average model, which is an ARMA$(0,q)$ model fitted over the time-series obtained by taking the logarithmic transformations of the original time-series observations. We refer readers to \citep{Box2015} for further details of ARMA.

\subsubsection{DeepAR}
Before discussing the details of how DeepAR can be used for sampling purposes, we briefly describe the building blocks of DeepAR, assuming familiarity with  
the basic concepts of neural networks, which can be found for instance in \citep{goodfellow2016deep}.

\paragraph{Recurrent neural network}
Due to the underlying limitation of multi-layer neural networks on processing sequential data, a specialized family of neural networks, namely recurrent neural networks (RNNs), has been developed. 
Let $\{d^{t}\}_{t \in \mathcal{T}}$ be a network input sequence and assume the data is fed in an ordinal manner. 
Given a vector $\red{\Theta}$ of parameters associated with the network and a state mapping function $g(\cdot)$, the state of the hidden units of the network at stage $t \in \mathcal{T}$, defined by
$h^{t} = g(h^{t-1}, d^{t}; \red{\Theta})$,
carries the encoded information of the previous steps by receiving $h^{t-1}$ as its input. The calculated states are then used for the output generation. 

\paragraph{DeepAR: probabilistic forecasting with autoregressive recurrent networks}
DeepAR is a neural network-based 
time-series forecasting method developed by \citet{Salinas2019}. It uses a specialized RNN architecture, called the long short-term memory, that is capable of carrying and forgetting information. It is a \emph{probabilistic} time-series forecasting method, i.e., its output defines the parameters of a probability distribution. Given history $\hat{\xipar}_{[t^*]}$, we assume $i \in \{1, \hdots, N\}$ is used to denote the sample index. In DeepAR, the neural network is trained to maximize the following objective function:
\begin{align*}
   \sum_{i = 1}^N \sum_{t = 0}^{t^*} \log (\ell(\hat{\xipar}_{i,t}|\theta(\mathbf{h}_{i,t}))
\end{align*}
where $\ell$ is the likelihood function, $\mathbf{h}_{i,t}$ represents network \red{hidden layer} parameters, and $\theta$ is a function that receives network \red{hidden layer parameters} an\red{d} returns parameters for the likelihood function.
Due to its construction, it can learn global patterns from many different time-series and covariates, where learned patterns are used to improve the forecasting performance for other time-series. For example, it can be used to predict weekly sales by inferring patterns from covariates such as weekly temperature records, boolean holiday indicators, among others, to improve the overall forecasting performance. \red{Besides DeepAR, there exists other recent deep learning-based probabilistic time series forecasting models (e.g., see Deep-state-space  \citep{rangapuram2018deep}). However, we only consider the combination of DeepAR and MSP policies, as DeepAR has been commonly used as a strong benchmark in many studies.
}

\paragraph{Conditional sample acquisition via ancestral sampling with DeepAR}
In a rolling\red{-}horizon framework, let $t^*$ denote the current stage and $\hat{\xipar}_{[t^*]}$ be the observed history. 
While making the decisions associated with stage $t^*$, a natural approach is to provide the optimization model \red{with} the updated forecasts for the future stages, i.e., the data obtained by retraining the neural network including the newly obtained observations, $\hat{\xipar}_{t^*}$. However, such an approach is not practical when training the network requires a significant computational effort, or the allowed time to get updated forecasts between decision-making stages is small. 
Although there exist many approaches to reduce training times, e.g., transfer learning, retraining might not always be a viable approach for some real-world applications. 

To cope with this problem, a heuristic sampling approach can be considered, by exploiting the ancestral sampling capacity of the encoder-decoder architecture of RNN. Ancestral sampling is a commonly used technique in the field of probabilistic graphical models \citep{goodfellow2016deep}. 
It sequentially generates samples from a complex joint distribution in such a way that a single sample at each step is obtained from a relatively easy-to-sample distribution compared to that of the joint density. Such a forecasting process of DeepAR is provided in Algorithm \ref{alg:sample_paths_DeepAR}, which returns $|\setSce|$ sample paths for the future, denoted by $\tilde{\xipar}^s, s \in \{1,\hdots,|\setSce|\}$. \red{Note that the algorithm takes as input the current decision-making stage, $t^*$, and the observed history of the uncertainty, $\hat{\xi}_{[t^*]}$, and generates sample paths without retraining the forecasting model. In other words, while using MSP combined with DeepAR in real-time decision making, future-stage predictions are automatically updated after each observation (i.e., after some uncertainty is revealed to the decision-maker) thanks to the time-step-dependent forecasting ability of the RNN-based DeepAR architecture.} 
\IncMargin{1em}
\begin{algorithm}
\caption{Conditional sample path generation with an RNN-based time-series forecasting method}
\label{alg:sample_paths_DeepAR}
\KwIn{$h^{t^*-1}, \red{\Theta}, \hat{\xipar}_{[t^*]}, \mathsf{D}(\cdot)$}  
\KwOut{$\{\tilde{\xipar}_t^s\}_{t \in \{t^*+1, \hdots, \Tpar\}, s \in \{1,\hdots,|\setSce|\}}$} 
\For{$s \in \{1, \hdots, |\setSce|\}$}{
    Let $\tilde{\xipar}^s_{t^*} = \hat{\xipar}_{t^*}$ \\
    \For{$t \in \{t^*, \hdots, \Tpar \red{-1} \}$}{
        $h^{t} = g(h^{t-1}, \red{\tilde{\xipar}_t^s}; \red{\Theta})$ \\
    	$\tilde{\xipar}^s_{\red{t+1}} = \textbf{Sample-from}(\mathsf{D}(h^{t}))$
    }
}
\end{algorithm}
\DecMargin{1em}

\noindent
In the algorithm, $\mathsf{D}(h)$ denotes any predefined probability distribution parametrized by $h$, which in our case is used to refer RNN output parameters.
We use $\textbf{Sample-from}(\cdot)$ to represent the function that returns a single draw from a given probability distribution. Although this function can be used to incorporate diverse sampling techniques, in this paper, we assume that it uses the crude Monte Carlo sampling.
\section{Experimental setup}\label{sec:evaluation}
As previously discussed, there exist many applications that require close coordination of time-series forecasting techniques and MSP. In this work, since our objective is to quantify the performance gap between the combination of general MSP with traditional and modern time-series forecasting techniques, we illustrate the results obtained on a particular MSP application, namely a multi-item lot-sizing problem with backlogging and production lag~\citep{daryalal2020lagrangian}, which contains both discrete and continuous variables. We conduct experiments using two well-known publicly available sales (demand) datasets. 

\subsection{Multi-item stochastic lot-sizing problem} 
The multi-item stochastic lot-sizing problem with backlogging and production lag (MSlag) aims to determine the optimal number of items to be produced at each time period (decision stage) while meeting customer demand for each item at each period, which is described as a stochastic process. Since the production cost at each period might differ and/or future demand may not be satisfied with the existing production capacity due to uncertainty, throughout the planning horizon, the decision maker is allowed to retain items in the inventory, as long as the inventory capacity restriction is satisfied, while spending holding costs for storing products in the inventory. Furthermore, in this problem setup, there exists a non-negative sequence describing fixed setup cost for deciding to produce any product at each period, regardless of the amount of the corresponding items to be produced \citep{chen2010applied}. In terms of time taken at production, there is a per period capacity, however overtime is allowed at the expense of incurring a higher production cost. Unmet demand at a period can be backlogged and satisfied at a future period. 

The MSP model for MSlag is provided below, and the detailed description for the parameters and decision variables are provided in Table \ref{t:MSlagNotation}.
\begin{itemize}[leftmargin=0.45cm]
    \item \textbf{Objective function. Minimize inventory, backlog, production setup, overtime costs}:
    \begin{align*}
        \hspace{-0.3cm}\min \ & \Exp\Bigg[\hspace{-0.1cm} \sum\limits_{t \in \Tset} \hspace{-0.1cm} \bigg( \hspace{-0.1cm} \sum\limits_{j \in \Jset} \hspace{-0.1cm} \Big(
C_{tj}^{i^+}(\xisupt)i^+_{tj}(\xisupt) \hspace{-0.07cm} + \hspace{-0.07cm} C_{tj}^{i^-}(\xisupt)i^-_{tj}(\xisupt) \hspace{-0.07cm} + \hspace{-0.07cm} C_{tj}^{y}(\xisupt)y_{tj}(\xisupt) \Big) \hspace{-0.07cm} + \hspace{-0.07cm}  C_{t}^{o}(\xisupt)o_{t}(\xisupt) \hspace{-0.1cm} \bigg) \hspace{-0.1cm}
\Bigg]
    \end{align*}
    \item \textbf{Constraints 1. State equations}: For all $t \in \Tset, j \in \Jset, \allxi$,
    \begin{align*}
        & i^-_{tj}(\xisupt)-i^+_{tj}(\xisupt) + i^+_{t-1,j}(\xisuptminusone) - i^-_{t-1,j}(\xisuptminusone) + x_{t-1,j}(\xisupt) = D_{tj}(\xisupt)
    \end{align*}
    \item \textbf{Constraints 2. Overtime measurements}: For all $t \in \Tset, \allxi$,
    \begin{align*}
        \sum\limits_{j \in \Jset}(TS_jy_{tj}(\xisupt) + TB_jx_{tj}(\xisupt)) - o_{t}(\xisupt)\leq C_t
    \end{align*}
    \item \textbf{Constraints 3. Link setup and production decisions}: For all $t \in \Tset, j \in \Jset, \allxi$,
    \begin{align*}
        M_{tj} \, y_{tj}(\xisupt) \geq x_{tj}(\xisupt)
    \end{align*}
    \item \textbf{Constraints 4. Capacity bounds}: For all $t \in \Tset, j \in \Jset, \allxi$,
    \begin{align*}
        & i^+_{tj}(\xisupt) \leq I_{tj} \\
        & i^+_{tj}(\xisupt) + x_{tj}(\xisupt) \leq I_{t+1,j} \\
        & o_{t}(\xisupt) \leq O_{t}
    \end{align*}
    \item \textbf{Constraints 5. Nature of variables}: For all $t \in \Tset, j \in \Jset, \allxi$,
    \begin{align*}
        & x_{tj}(\xisupt),i^+_{tj}(\xisupt), i^-_{tj}(\xisupt), o_{t}(\xisupt) \geq 0 \\
        & y_{tj}(\xisupt) \in\{0,1\}
    \end{align*}
\end{itemize}

\noindent Note that in this optimization model, we assume that there is no production related cost. More detailed explanations of the model can be found in \citep{daryalal2020lagrangian}.

\begin{table}[!ht]
\caption{Notation used in the MSlag model}
\label{t:MSlagNotation}
\small
\begin{tabular}{ll ll} 
\toprule
\multicolumn{4}{l}{\textbf{\textsc{Sets and Indices:}}} \\
$\Tset$ & Time periods, $t \in \Tset$ \\
$\Jset$ & Product types, $j \in \Jset$ \\*[0.05cm] 
\midrule
\midrule
\multicolumn{4}{l}{\textbf{\textsc{Parameters:}}} \\
\multicolumn{2}{l}{\it Deterministic:} & \multicolumn{2}{l}{\it Uncertain:} \\ 
$M_{tj}$: & Big-M value for modeling & $\xisupt$ & Vector of all random variables up to period $t$ \\
$TS_j$ & Setup time & $D_{tj}(\xisupt)$ & Demand of products \\
$TB_j$: & Unit production time & $C_{tj}^{y}(\xisupt)$ & Production cost \\
$C_t$: & Production capacity & $C_{tj}^{i^-}(\xisupt)$ & Backlog cost \\
$I_{tj}$: & Inventory capacity & $C_{tj}^{i^+}(\xisupt)$ & Inventory cost \\
$O_{t}$: & Overtime bound & $C_{t}^{o}(\xisupt)$ & Overtime cost 
\\*[0.05cm] 
\midrule
\midrule
\multicolumn{2}{l}{\textbf{\textsc{Decision variables:}}} \\*[0.1cm]
$x_{tj}(\xisupt)$ & Production level \\
$i^+_{tj}(\xisupt)$ & Inventory level \\
$i^-_{tj}(\xisupt)$ & Backlog level \\
$o_t(\xisupt)$ & Overtime measurement \\
$y_{tj}(\xisupt)$ & \multicolumn{3}{l}{Production  decision which takes value $1$ if production is setup, 0 otherwise}
\\*[0.15cm]
\bottomrule
\hline
\end{tabular}
\end{table}

\subsection{Data description}
Regarding the stochastic process $\{D_{tj}(\xisupt)\}_{t \in \Tset, j \in \Jset}$, for demand scenario generation purposes, we used two publicly available datasets: Walmart Sales Forecasting data~\citep{Kaggle_Walmart} and Corporaci\'on Favorita Grocery Sales Forecasting data~\citep{Kaggle_Favorita}. For the remaining optimization model parameter generation, we followed the steps described in \citep{daryalal2020lagrangian} which builds on the steps provided in \citep{helber2013dynamic}. We next briefly outline the characteristics of the used demand data. 

The Walmart dataset contains 12 columns of weekly sales data records for each department at each store, from 2010-02-05 to 2012-11-01. We let each combination of ``\textit{(store, department)}" be one demand time-series. Given the dataset, we decided to drop seven columns due to a large amount of missing data and a low correlation with the other features. We note that removal of these features would not lead to any predictive performance degradation, and producing robust and reliable statistical forecasting results using missing data is beyond the scope of this paper. We directly used and/or modified the remaining features to obtain new covariates to be used as an input for DeepAR: Year, Month, Week, Day, Temperature, Fuel price, and Holiday indicator. Despite the fact that temperature and fuel price might not be the most suitable covariates as they include data on future events, for this research, we assumed both time series as given, thanks to the advances in the field of time-series forecasting. For the sake of prediction, we extracted the data of the last eight weeks and trained all prediction models for the data without the last eight weeks.  For (store, department) combinations to be used in our numerical evaluation, we randomly sampled from a subset of the original dataset where the forecasting methods had a noticeable predictive performance difference.

\red{For the Walmart and Favorita datasets, the planning horizons correspond to the time steps $\{135, \hdots, 142\}$ and $\{580, \hdots, 591\}$, as can be observed in Figure \ref{fig:forecasting_walmart} and Figure \ref{fig:forecast_favorita}, respectively.} We use all the data available before the beginning of the planning horizon in the training phase. In the rest of paper, we use $\Tpred$ to denote the \emph{prediction horizon length} for DeepAR.

\subsection{Implementation details}
All the MSP model policy generation and evaluation algorithms are implemented in Python, and Gurobi 9.0.0. is used as a mixed-integer programming solver. Experiments are conducted on a MacOS workstation with Quad-Core 3GHz Intel i5-8500B CPU and 16 GB memory. 


Regarding DeepAR, a PyTorch \citep{NEURIPS2019_9015} (a deep learning package for Python) implementation developed by \citet{zhang2019} was used as the baseline code. We note that there exists an open source package called GluonTS~\citep{gluonts_jmlr}, which is capable of performing deep learning-based time-series forecasting in an integrated and automated manner. However, to the best of our knowledge, the conditional sample path generation feature is not supported by GluonTS.

\section{Numerical results}\label{sec:experiments}
In this section, we first compare the predictive performance of the three forecasting methods. In the policy evaluation section, we analyze the impact of improved forecasts on the deterministic and two-stage policies. For the former, in addition to the risk-neutral setting, we consider a risk-averse approach \red{where the demand values are replaced by the highest value of their corresponding prediction interval.}


\subsection{Forecasting performance}\label{sec:forecastingPerformance}
We evaluate the overall predictive performance of AR(1), moving average and DeepAR models on the used datasets. Specifically, for the Walmart dataset and the Favorita dataset, we use regular moving average (MA) and logaritmic moving average (LogMA) models, where we consider the moving average parameters of 8 and 12, respectively\red{, both obtained via grid search.} Our choice of LogMA rather than MA for the Favorita dataset is based on the fact that it has been empirically shown to be competitive with the state-of-the-art methods~\citep{Kaggle_Favorita}. In terms of the effort required to build the forecasting models, the parameters associated with ARIMA-based models are calculated in less than ten minutes, while training DeepAR is significantly more time consuming, requiring around twelve hours.

\red{
In order to train the DeepAR model, we split the dataset into three subsets: training, validation, and test. 
For the Walmart dataset, training, validation, and test datasets involve time steps $\{0, \hdots, 126\}$, $\{127, \hdots, 134\}$, and $\{135, \hdots, 142\}$, respectively. 
Similarly, for the Favorita dataset, training, validation, and test datasets involve time steps $\{0, \hdots, 529\}$, $\{530, \hdots, 579\}$, and $\{580, \hdots, 591\}$, respectively. 
We performed hyper-parameter tuning for the model parameters, where the model is trained with the training dataset, and its performance is measured based on the validation dataset.
The obtained DeepAR hyper-parameters are provided in Table~\ref{tab:deeparhp}.
Using the identified parameters, we performed time series forecasting for the time steps corresponding to the test dataset.
}

\begin{table}[ht!]
\cred
\small
\centering
\caption{\cred DeepAR hyper-parameters}
\label{tab:deeparhp}
\begin{tabular}{l|cccc}
          & \multicolumn{1}{l}{\textbf{Batch size}} & \multicolumn{1}{l}{\textbf{Learning rate}} & \multicolumn{1}{l}{\textbf{LSTM layers}} & \multicolumn{1}{l}{\textbf{LSTM units}} \\ \hline
\textbf{Favorita} & 64                                      & 1e-3                                       & 3                                        & 40                                      \\
\textbf{Walmart}  & 64                                      & 1e-3                                       & 6                                        & 80                                     
\end{tabular}
\end{table}

\begin{figure}[!ht]
\centering
\includegraphics[width=1.0\textwidth]{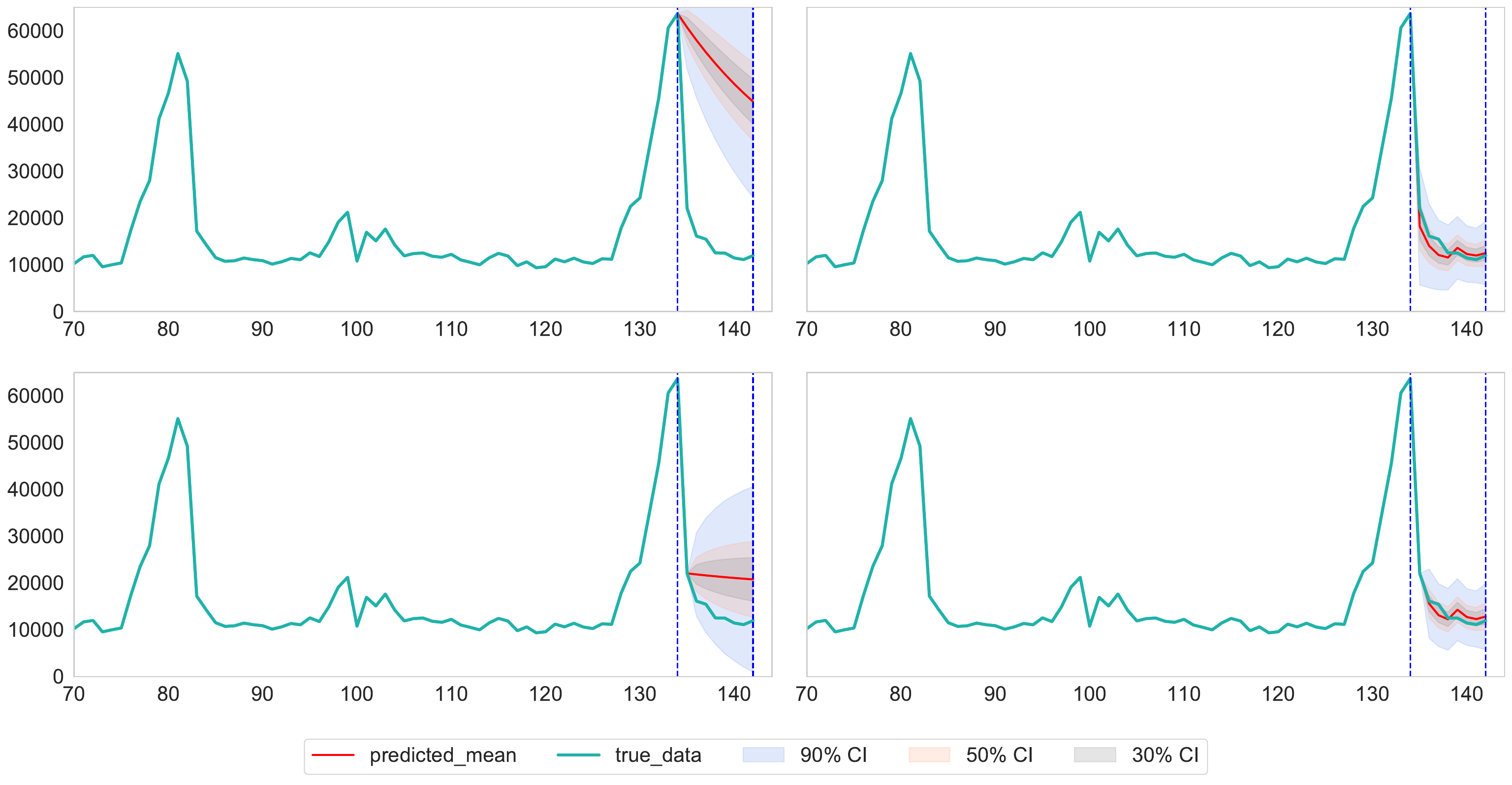}
\caption{Probabilistic product demand forecasting for a selected time-series from the Walmart data (top-left: unconditioned AR(1), bottom-left: 1-step conditioned AR(1), top-right: unconditioned DeepAR, bottom-right: 1-step conditioned DeepAR). \red{The first decision-making stage, $t = 0$, corresponds to the time step 135.} Dotted vertical blue lines indicate the beginning and the end of the planning horizon for MSlag.}
\label{fig:forecasting_walmart}
\end{figure}

We use a rolling\red{-}horizon framework to evaluate different policies, and thus adopt conditional sampling. In rolling\red{-}horizon framework, at every decision-making stage, that is, for any MSP model to be solved, the forecasts are obtained using available history, i.e., all the true observations made up to that stage. For the first stage (i.e., for the first MSP to be solved), forecasts are called \emph{unconditioned} (or \emph{0-step conditioned}); for the second one they are referred to as \emph{1-step conditioned}, and so on. For $t'$-step conditioned case, we refer to $t'$ as the \emph{conditioning step}. In order to provide a better idea on the notion of conditioning, we illustrate unconditioned and one-step conditioned probabilistic forecasts of AR(1) and DeepAR for one particular time series from the Walmart dataset in Figure~\ref{fig:forecasting_walmart}.

Figure~\ref{fig:forecasting_walmart} demonstrates the poor-quality performance of AR(1) compared to that of DeepAR, for both unconditioned and 1-step conditioned forecasts, in terms of both the predicted mean values and confidence interval widths. Although we provide the forecasting results for only one particular product in this figure as an example, we note that the plots for the remaining products look very similar. In what follows, we present a formal overall performance measure and compare all three considered methods.

Since DeepAR is used as one of our main forecasting methods, we rely on the so-called normalized deviation (ND) performance metric~\citep{Salinas2019}. However, as predictions are performed for different step conditions in our rolling-horizon framework, namely from unconditioned to $(\Tpar-1)$-step conditioned if the planning horizon length is $\Tpar$, we define and use ND values for each possible conditioning step. 
More specifically, given $\Tpar$ as the length of the planning horizon, and a conditioning step $t' \in \{0,\hdots, \Tpar-1 \}$, the \emph{$t'$-step ND value} is computed as
\begin{align*}
\text{ND$(t',\Tpar)$} := \frac{\sum_{t = t'+1}^\Tpar \sum_{j \in \Jset} |D_{tj} - \hat{D}_{tj}|}{\sum_{t = t'+1}^\Tpar \sum_{j \in \Jset} |D_{tj}|}
\end{align*}
where for (the time-series corresponding to) product $j$ at time $t$, $D_{tj}$ and $\hat{D}_{tj}$ denote the true and predicted median demand values.

We report the obtained ND values of selected set of time series, for each conditioning time step for the Walmart and Favorita datasets in Table~\ref{tab:nd_results}, where column labels refer to the conditioning steps ($t'$). 
\setlength{\tabcolsep}{4pt} 
\renewcommand{\arraystretch}{1.13} 
\begin{table}[!htp]
\small
\centering
\caption{{ND($t',\Tpar$) values at each prediction period for different forecasting methods}}
\label{tab:nd_results}
    \subfloat[{Walmart dataset, $\Tpar = 8$, $\Tpred = 8$} \label{tab:walmart_nd}]{
\begin{tabular}{lrrrrrrrr}
\toprule
\textbf{Method}\ /\ $t'$ &     \textbf{0} &     \textbf{1} &     \textbf{2} &     \textbf{3} &     \textbf{4} &     \textbf{5} &     \textbf{6} &     \textbf{7} \\
\midrule
\textbf{AR(1)}  & 0.34 &  0.36 &  0.28 &  0.18 &  0.16 &  0.15 &  0.17 &  0.14 \\
\textbf{MA}  &  0.25 &  0.25 &  0.19 &  0.14 &  0.13 &  0.12 &  0.13 &  0.13 \\
\textbf{DeepAR} &  0.17 &  0.14 &  0.13 &  0.12 &  0.13 &  0.12 &  0.14 &  0.15 \\
\bottomrule
\end{tabular}
} \\[0.6em]
    \subfloat[{Favorita dataset, $\Tpar = 12$, $\Tpred = 12$} \label{tab:favorita_nd}]{
\begin{tabular}{lrrrrrrrrrrrr}
\toprule
\textbf{Method}\ /\ $t'$ &     \textbf{0} &     \textbf{1} &     \textbf{2} &     \textbf{3} &     \textbf{4} &     \textbf{5} &     \textbf{6} &     \textbf{7} &     \textbf{8} &     \textbf{9} &     \textbf{10} &     \textbf{11} \\
\midrule
\textbf{AR(1)}  &  0.99 &  1.22 &  0.56 &  0.56 &  0.60 &  0.60 &  0.56 &  0.55 &  0.50 &  0.41 &  0.37 &  0.31 \\
\textbf{LogMA}  &  0.50 &  0.47 &  0.44 &  0.43 &  0.44 &  0.43 &  0.42 &  0.42 &  0.41 &  0.37 &  0.35 &  0.31 \\
\textbf{DeepAR} &  0.48 &  0.46 &  0.43 &  0.39 &  0.40 &  0.40 &  0.41 &  0.42 &  0.43 &  0.41 &  0.37 &  0.45 \\
\bottomrule
\end{tabular}
}
\end{table} 

Table~\ref{tab:walmart_nd} shows that \red{for the Walmart dataset} DeepAR performs the best \red{for the first four time steps, while it is competitive with the other methods for the last four time steps, namely $t' \in \{4,5,6,7\}$. In contrast,}  AR(1) \red{performs} the worst for the Walmart dataset \red{for most of the time steps}. For the Favorita dataset, as seen in Table \ref{tab:favorita_nd}, while DeepAR outperforms the others \red{for the majority of the time steps}, the performance gap between LogMA and DeepAR is not as large as in the Walmart case, which is inline with the reported performance of LogMA in \citep{Kaggle_Favorita}. 
Note that the forecasting performance of all the methods for stages closer to the end of the planning horizon, i.e., the ones obtained by using a large conditioning step value, are highly accurate. 
This can be attributed to the fact that the value of the target time-series stabilizes as time index increases, while there is some unstable behavior early on, e.g., a huge peak is observed at the beginning of the horizon (see Figure \ref{fig:forecasting_walmart}). It is also important to recognize that time-series forecasting methods usually perform well for short-term prediction, while they frequently fail to provide high-accuracy long-term predictions. 

As the ND values are computed based on predicted median demand values, they provide insights mostly on the quality of point forecasts, which play an important role in risk-neutral deterministic policies. On the other hand, for two-stage policies, as well as for risk-averse deterministic policies, we instead use probabilistic forecasts. That is, we rely on the provided confidence intervals to either sample scenarios or estimate the worst-case scenario/build an uncertainty set. Therefore, we use an additional metric to be able to assess the overall performance of the considered methods in terms of their probabilistic forecasts. In that regard, we rely on the $\rho$-risk (quantile loss) metric \citep{Salinas2019}, which is used to quantify the forecasting accuracy of a quantile $\rho$ of the predictive probability distribution.

Given $\rho \in (0,1)$, the $\rho$-risk of a product $j \in \Jset$ is defined as 
\begin{align*}
L_{\rho} (Z_j, \hat{Z}^{\rho}_j) = 2 \cdot (\hat{Z}^{\rho}_j - Z_j) \cdot \left(\rho \cdot \mathbb{I}_{(\hat{Z}^{\rho}_j > Z_j)} - (1-\rho) \cdot \mathbb{I}_{(\hat{Z}^{\rho}_j \leq Z_j)} \right)
\end{align*}
where $\mathbb{I}_{(\cdot)}$ is the indicator function. $Z_j = \sum_{t_0+L}^{t_0 + L + S} D_{tj}$ where $t_0$ is the forecast start point, $L$ is the lead time after the forecast start point, and $S$ is the length of the prediction range. In our case, $t_0$ corresponds to the start of the planning horizon (i.e., the period index of our $t=0$), there is no lead time, i.e., $L= 0$, and the length of the prediction range is equal to the length of the planning horizon, i.e., $S = \Tpar$. Lastly, the predicted $\rho$-quantile value of $Z_j$ is denoted by $\hat{Z}_j^{\rho}$. Note that having a lower value of $\rho$-risk is better. Analogously, $\rho$-risk for the entire set of products, $\rho$-risk, is computed by using the following formula:
\begin{align*}
\frac{\sum_{j \in \Jset} L_{\rho}(Z_j, \hat{Z}^{\rho}_j)}{\sum_{j \in \Jset} Z_j}
\end{align*}

For our datasets, using a high confidence level of $\rho = 0.9$, for the Walmart dataset, we obtain the $\rho$-risk values of 1.57, 1.57 and 1.08 for AR(1), MA and DeepAR, respectively. On the other hand, for the Favorita dataset, the 0.9-risk values are attained as 2.32, 1.65 and 1.27 respectively for AR(1), MA and DeepAR. \red{These values indicate that DeepAR is able to generate more accurate probabilistic forecasts with lower variance.}

Lastly, we note that we are not aware of any previous studies reporting the forecasting performance values (e.g. $\rho$-risk) of the Walmart and Favorita datasets. However, the above-mentioned values seem to be inline with the ones provided in Table 1 of \citep{Salinas2019}. For a challenging natural number valued dataset, namely \textsc{parts}, the authors report 0.9-risk value as 1.06 when $L=0$ and $S=8$ (as in our Walmart dataset). Also, for another dataset, called \textsc{electricity}, which can be considered as a favorable dataset for forecasting due to inherent seasonality, they find the 0.9-risk as 1.33 when $L=3$ and $S=12$ (as in our Favorita set). 


\subsection{Policy evaluation}

In this section, we present the quality of different policies combined with different forecasting methodologies. More specifically, we focus on deterministic and two-stage policies (see Section \ref{subsec:det} and Section \ref{subsec:stoch}) as well as a risk-averse (robust) policy optimizing the worst-case cost instead of the expected cost.

As a key performance measure we consider the percentage gap between the perfect information (PI) bound and the objective function value obtained by evaluating considered policies in a rolling-horizon framework. We refer to this value as ``Gap \%". The PI bound can be obtained by solving the optimization problem described in Section \ref{subsec:det}, by replacing $\xipar^{\text{expected}}_{t}$ with $\hat{\xipar}_{t}$ (i.e., replacing the expected value with the true observation) for every stage $t \in \Tset$. The choice of the PI bound as a reference point is common in the MSP literature, and stems from the facts that the optimal solution to an MSP is not possible to obtain in general, that the PI bound is computationally quite cheap to compute, and that finding improved bounds are usually computationally demanding. In our experiments, we consider multiple product groups, where a group corresponds to a particular set $\Jset$ in MSlag, as such we report the Gap \% values averaged over all the product groups.

\red{We note that the MSP model has relatively complete recourse. As such inaccurate forecasts do not lead to any infeasibility. However, in such cases, backlog penalties are incurred, which increases the model objective function values, implying larger Gap \% values.}

\subsubsection{Deterministic risk-neutral setting}

\begin{figure}[h]
    \centering
    \subfloat[the Walmart dataset \label{fig:Walmart1}]{
    \includegraphics[width=0.5\textwidth]{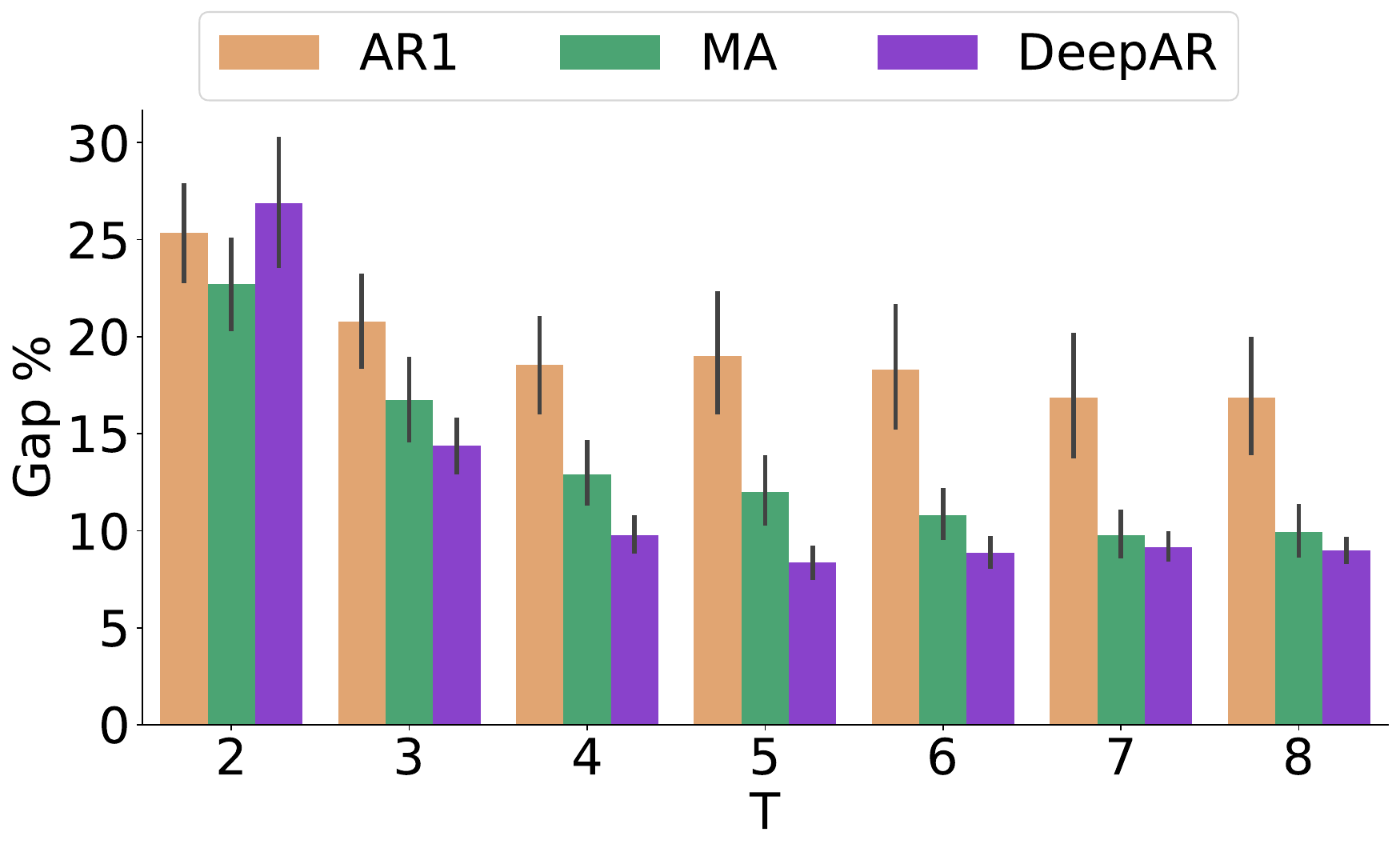}
    }
    \hfill
    \subfloat[the Favorita dataset \label{fig:Favorita1}]{
    \includegraphics[width=0.5\textwidth]{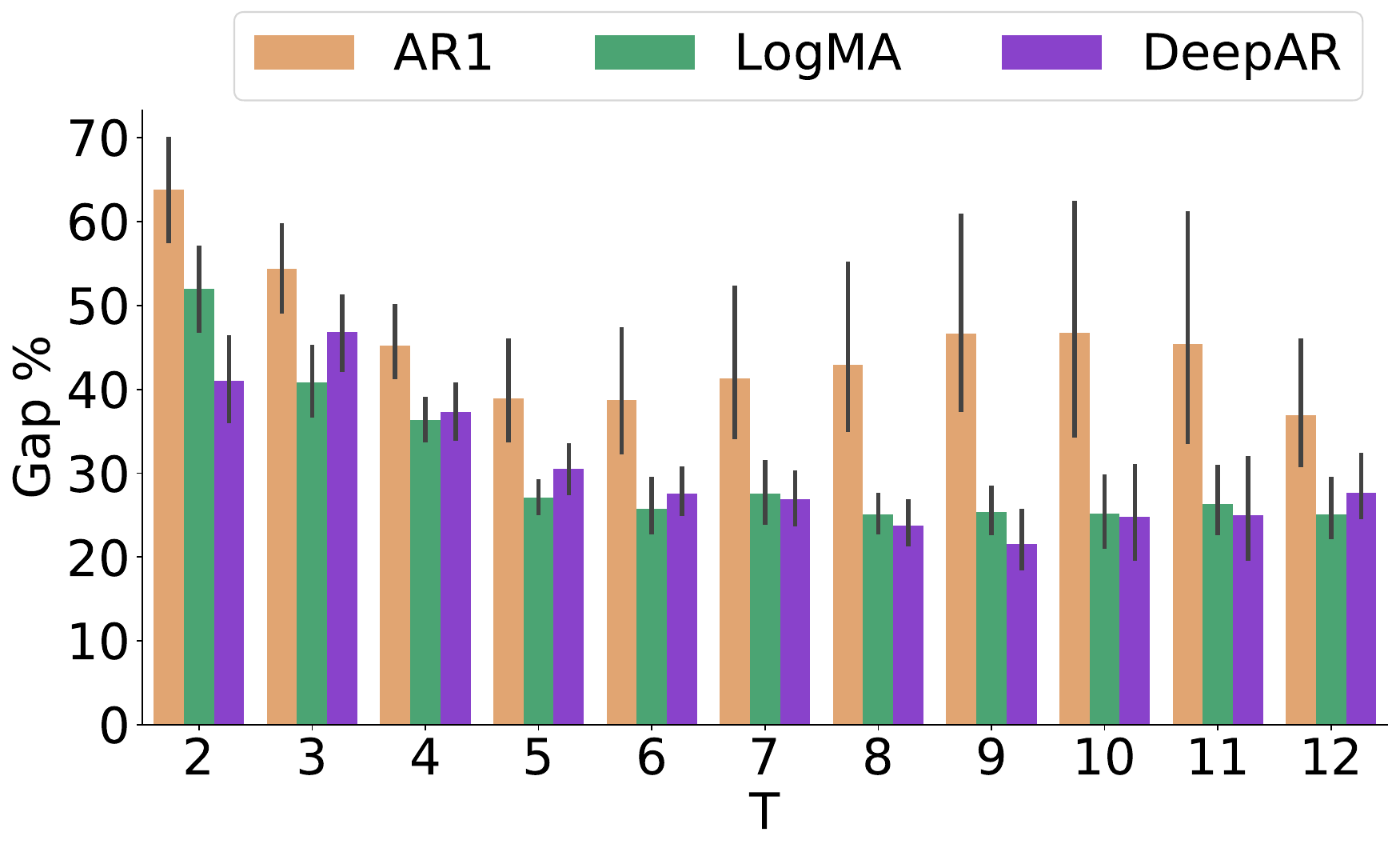}
    }
    \caption{Gap \% values for different planning horizons ($T$) under risk-neutral setting}
    \label{fig:deterministic_gapplots}
\end{figure}

In the rolling-horizon experiments with the deterministic policy, we consider 50 product groups (i.e., $\Jset$ sets) and use the following parameter settings:
\begin{itemize}
    \item The Walmart dataset: $|\Jset| = 10$ and $\Tpar \in \{2,3, \hdots, 8\}$
    \item The Favorita dataset: $|\Jset| = 15$ and $\Tpar \in \{2,3, \hdots, 12\}$
\end{itemize}
We provide the results obtained with AR(1), (Log)MA, and DeepAR forecasting methods for varying length of the planning horizon in Figure~\ref{fig:deterministic_gapplots}, which also includes error bars for each case, showing the deviation in different product groups.

We observe that, for the deterministic policy, overall, DeepAR leads to the lowest Gap \% values, while AR(1) consistently performs the worst. We also note that the benefits of improved forecasting can be substantial. For instance, for the Walmart dataset with $\Tpar = 8$, the average Gap \% values for AR(1), MA and DeepAR are 16.8\%, 10.0\% and 9.0\%, respectively. Furthermore, for the same settings (Walmart with $\Tpar = 8$), average of the worst and best three product group Gap \% values for AR(1) are 46.6\% and 7.0\%, whereas, DeepAR leads to 12.3\% and 6.5\% gaps for the same product groups, which shows that a significantly higher-quality deterministic policies can be obtained with the help of improved forecasts. For both datasets, we find that MA-based methods and DeepAR perform similarly, which is consistent with their predictive performance reported in Section \ref{sec:forecastingPerformance}. In addition, these two methods show low variability in Gap \% values, and, in general, their average performance as well as the variance improve as $\Tpar$ gets closer to the prediction horizon length, $\Tpred$. We note that our analysis with varying number of products ($|\Jset|$) for fixed $\Tpar$ values did not yield substantially different insights, thus was omitted for the sake of brevity.

\begin{figure}[h]
    \centering
    \subfloat[ \label{fig:spw1}]{
    \includegraphics[width=0.8\textwidth]{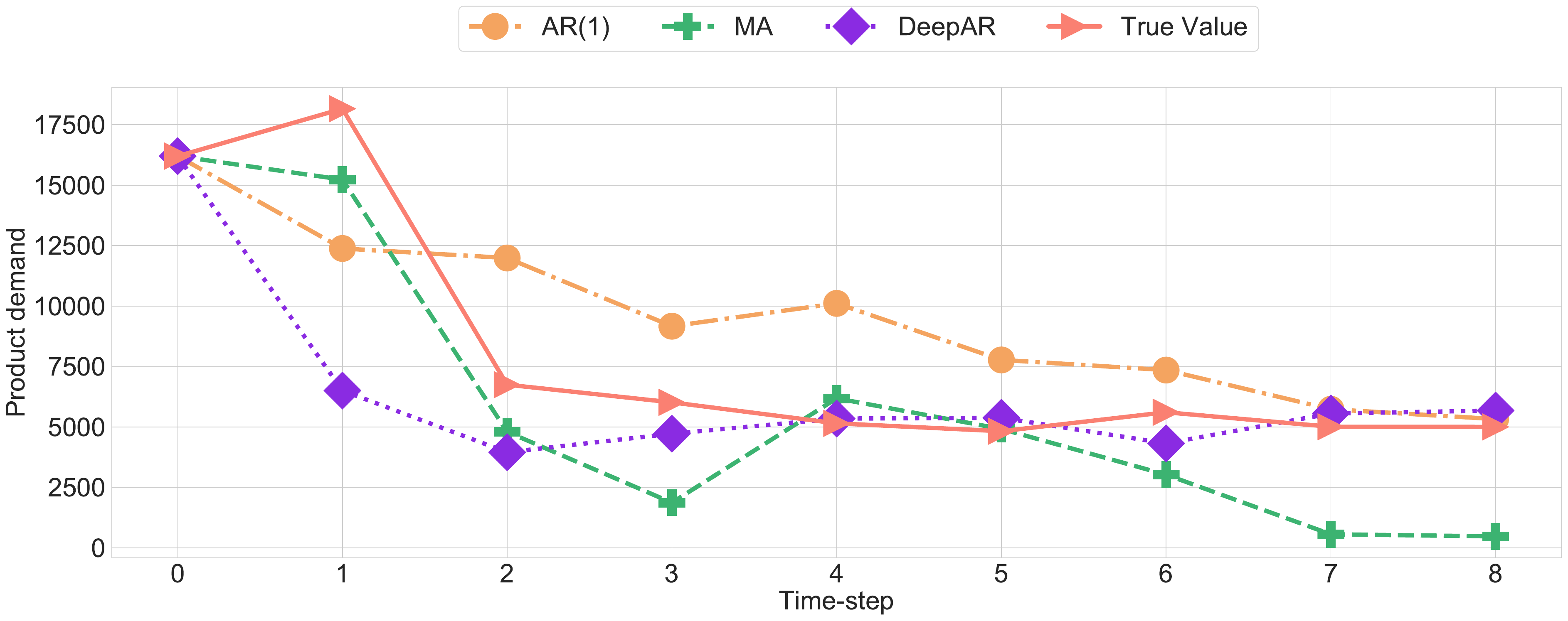}
    } \\
    \subfloat[ \label{fig:spw2}]{
    \includegraphics[width=0.8\textwidth]{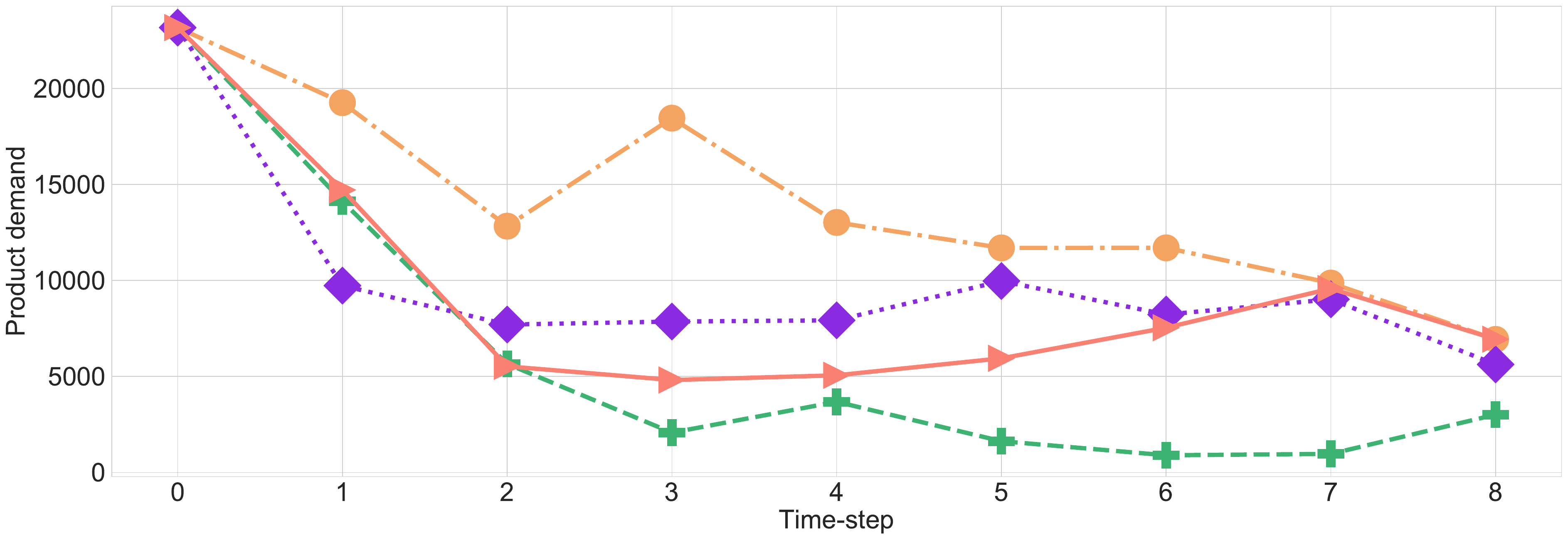}
    }
    \caption{Forecasting results for two different time-series from the Walmart dataset.}
    \label{fig:sample_paths_walmart}
\end{figure}

We observe a poor performance of DeepAR for the Walmart dataset when $\Tpar=2$. To elaborate, we provide Figure \ref{fig:sample_paths_walmart} which illustrates the forecasting result comparison of a sample path generated from each forecasting method, for two different time-series.
In Figure \ref{fig:spw1}, after the large peak at $t = 0$ (corresponding to time index 135 in the original time-series) as shown in Figure \ref{fig:forecasting_walmart}, DeepAR predictions tend to rapidly decrease with respect to the true values, producing some errors in the first few steps of the forecasting horizon. On the other hand, after the first few steps, DeepAR produces significantly better results compared to AR(1) and MA. These forecasting results help explaining why DeepAR performs worse when $T=2$, while performing well overall. We also note that a similar phenomenon is observed in Figure \ref{fig:spw2} and in many other time-series, where in this case, DeepAR failed to perform better than MA. 

We also provide the average performance values for the Walmart dataset when $\Tpar = 2$ in Table \ref{tab:wal_T2}, which points to the inferior performance of DeepAR for $t'\in \{0,1\}$. Similar to Table \ref{tab:nd_results}, column indices of Table \ref{tab:wal_T2} correspond to the time step where the forecasting method was conditioned. 
\begin{table}[ht!]
\small
    \centering
    \caption{ND($t',\Tpar$) values for the Walmart dataset when $\Tpar = 2$.}
    \label{tab:wal_T2}
    \begin{tabular}{lrr}
    \toprule
    {Method} &  $t'=0$ &  $t'=1$ \\
    \midrule
    AR(1)  &  0.22 &  0.23 \\
    MA  &  0.20 &  0.22 \\
    DeepAR &  0.26 &  0.19 \\
    \bottomrule
    \end{tabular}
\end{table}

In the Appendix, we provide the ND($t',\Tpar$) values for all $\Tpar$ and $t' < \Tpar$ combinations, of the DeepAR model used in our analysis (i.e., the one trained with $\Tpred$ equal to the largest $\Tpar$ option) for both datasets. We note that for $t'=0$, the only case where DeepAR performance was worse than the other methods was for the Walmart case with $\Tpar=2$. 

\subsubsection{Deterministic risk-averse setting}
In this section, we present the deterministic policy evaluation results under a risk-averse setting. Specifically, for each forecasting method, we generate a $90\%$ \red{prediction} interval for the predictions of each time series, and use its maximum value to be taken as the value of the demand in the worst case. These worst-case demand estimates are provided to the MSlag models to be solved in the rolling\red{-}horizon framework. 

\begin{figure}[!htp]
    \centering
    \subfloat[the Walmart dataset \label{fig:Walmart2}]{
    \includegraphics[width=0.5\textwidth]{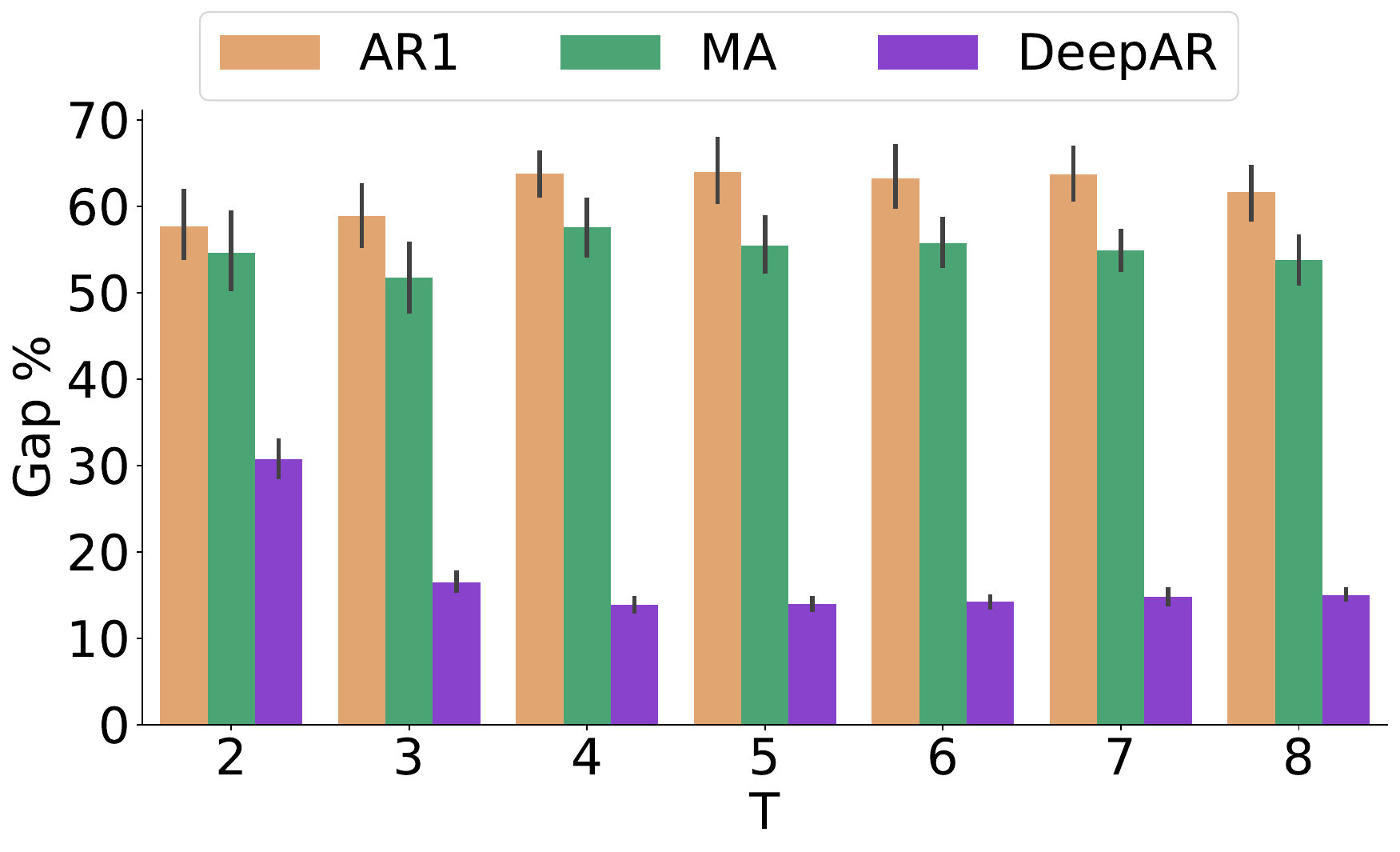}
    }
    \hfill
    \subfloat[the Favorita dataset \label{fig:Favorita2}]{
    \includegraphics[width=0.5\textwidth]{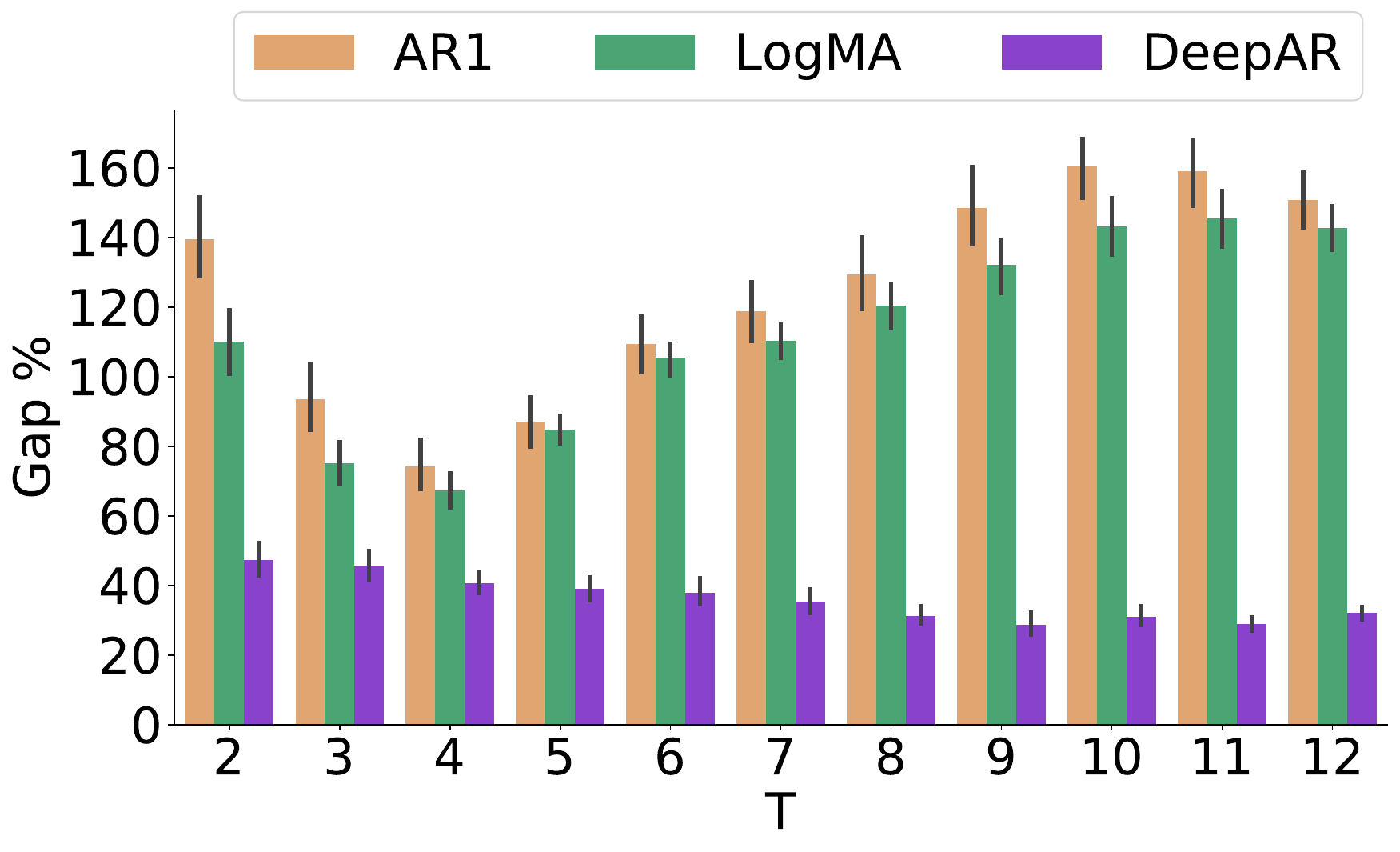}
    }
    \caption{Gap \% values for different planning horizons ($T$) under risk-averse setting}
    \label{fig:deterministic_riskaverse_gapplots}
\end{figure}

Figure \ref{fig:deterministic_riskaverse_gapplots} shows the obtained results with the risk-averse setting for varying planning horizon lengths. 
We observe that the policies generated with DeepAR significantly outperform the others. In particular, for the Favorita dataset, while LogMA-based deterministic policies are comparable to the DeepAR-based ones in the risk-neutral setting, LogMA loses its strength in the risk-averse setting. For instance, for the Favorita dataset with $\Tpar = 12$, the average Gap \% values are 150.7\%, 142.7\%, and 32.7\% for AR(1), LogMA, and DeepAR models, respectively.

\begin{figure}[htp!]
\centering
\includegraphics[width=1\textwidth]{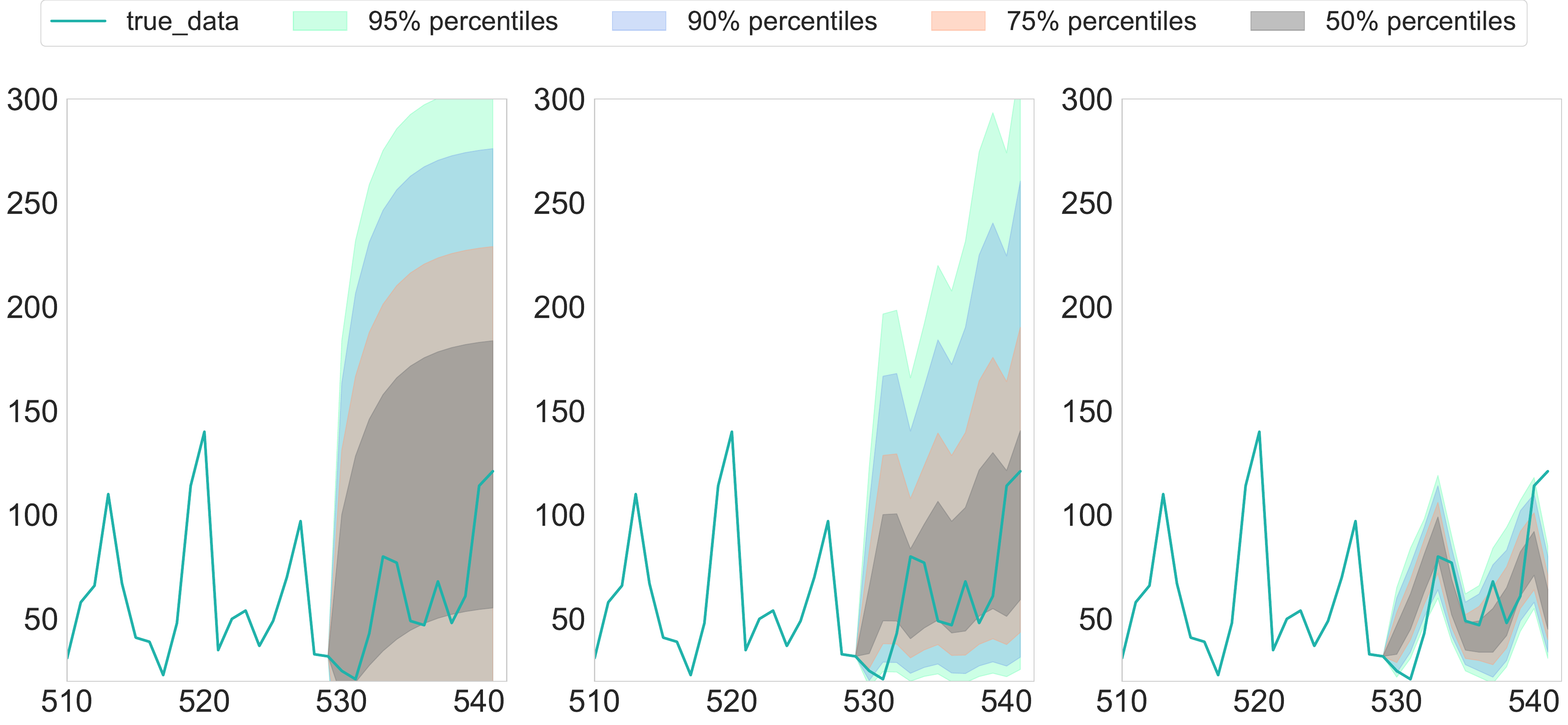}
\caption{A graphical representation of forecasting results of AR(1), LogMA and DeepAR for \red{the first portion of validation dataset of the} Favorita dataset, from left to right respectively, where its x-axis is time steps and y-axis is product demand}
\label{fig:forecast_favorita}
\end{figure}

As discussed by \citet{Salinas2019}, DeepAR tends to yield tighter prediction intervals. Such a phenomenon can also be observed from, for instance, Figure~\ref{fig:forecast_favorita} which presents sample forecasting results from the Favorita dataset, showing the significant reduced-variance forecasting power of DeepAR. We note that the $\rho$-risk measure analysis provided in Section~\ref{sec:forecastingPerformance} supports this finding.
Accordingly, the worst-case predictions of DeepAR are typically closer to the ground truth, which considerably benefits deterministic robust policies. 
In this regard, we note that considering a sophisticated forecasting method can be highly beneficial not only in terms of uncertainty modeling for stochastic programming, but also to construct smaller and accurate uncertainty sets for robust optimization models.

Lastly, regarding the results presented in Figure \ref{fig:deterministic_gapplots} and Figure \ref{fig:deterministic_riskaverse_gapplots}, note that the risk-averse policies typically yield larger Gap \% values than the risk-neutral ones since they are more conservative. However, in our experiments with DeepAR, they turn out to be quite close to each other. For example, for the Walmart dataset with $\Tpar = 8$, Gap \% values are 9.0\% and 15.0\% for risk-neutral and risk-averse settings, respectively, while for the Favorita dataset with $\Tpar = 12$, these values are 27.7\% and 32.1\%.

\subsubsection{Two-stage setting}

Finally, we present our findings on the impact of improved forecasts for two-stage policies using a rolling-horizon approach. In these experiments, we consider the number of scenarios, $|\setSce|$, as 9 and 15 for the Walmart and Favorita datasets, respectively. Unlike the deterministic policy experiments, for both datasets, we consider 30 product groups with $|\Jset| = 3$ and a range of planning horizon lengths $\Tpar \in \{2,3,4,5\}$. \red{Such a restricted setup consideration is mainly because two-stage stochastic models to be solved in the rolling horizon framework are significantly more computationally demanding than the deterministic models. We note that as we have mixed-binary recourse decisions, we solve the two-stage stochastic programs via their extensive form; developing a problem-specific solution method (e.g., decomposition algorithm) is out of scope of the paper.} 

\red{We chose the number of scenarios (9 and 15 for the two datasets) via preliminary experiments. We tested different number of sample sizes by enforcing a time limit of three hours. The provided sample sizes are the maximum ones that could be handled given the time limit, which is for all the product groups and the full rolling-horizon framework, rather than for solving a single two-stage stochastic program. More specifically, in our rolling-horizon experiments for each of the 30 product groups, we simulate a sample path by solving $T-1$ two-stage stochastic programs. As such, although we can solve a single two-stage stochastic program with more than 9 and 15 scenarios in the Walmart and Favorita cases, respectively, considering the total computational effort of solving $30(T-1)$ stochastic programs, we opted for a smaller sample size. It is also worth to mention that the low-variance forecasts provided by DeepAR would reduce the need for a large number of scenarios to get high-quality solutions for the two-stage stochastic programs.  
}

A particularly useful feature of DeepAR for our purposes is that it can generate sample paths with all non-negative demand values by appropriately adjusting the distribution assumption of the neural network architecture (e.g., using negative binomial or Poisson distribution for natural number time-series predictions). This is not the case in ARMA models, since the white noise term is not bounded. Therefore, as product demand is assumed to be non-negative valued, for sample paths obtained by the ARMA-based models, we replace negative prediction values with zeros.

\begin{figure}[htp!]
    \centering
    \begin{subfigure}[b]{0.48\textwidth}
        \centering
        \includegraphics[height=1.8in]{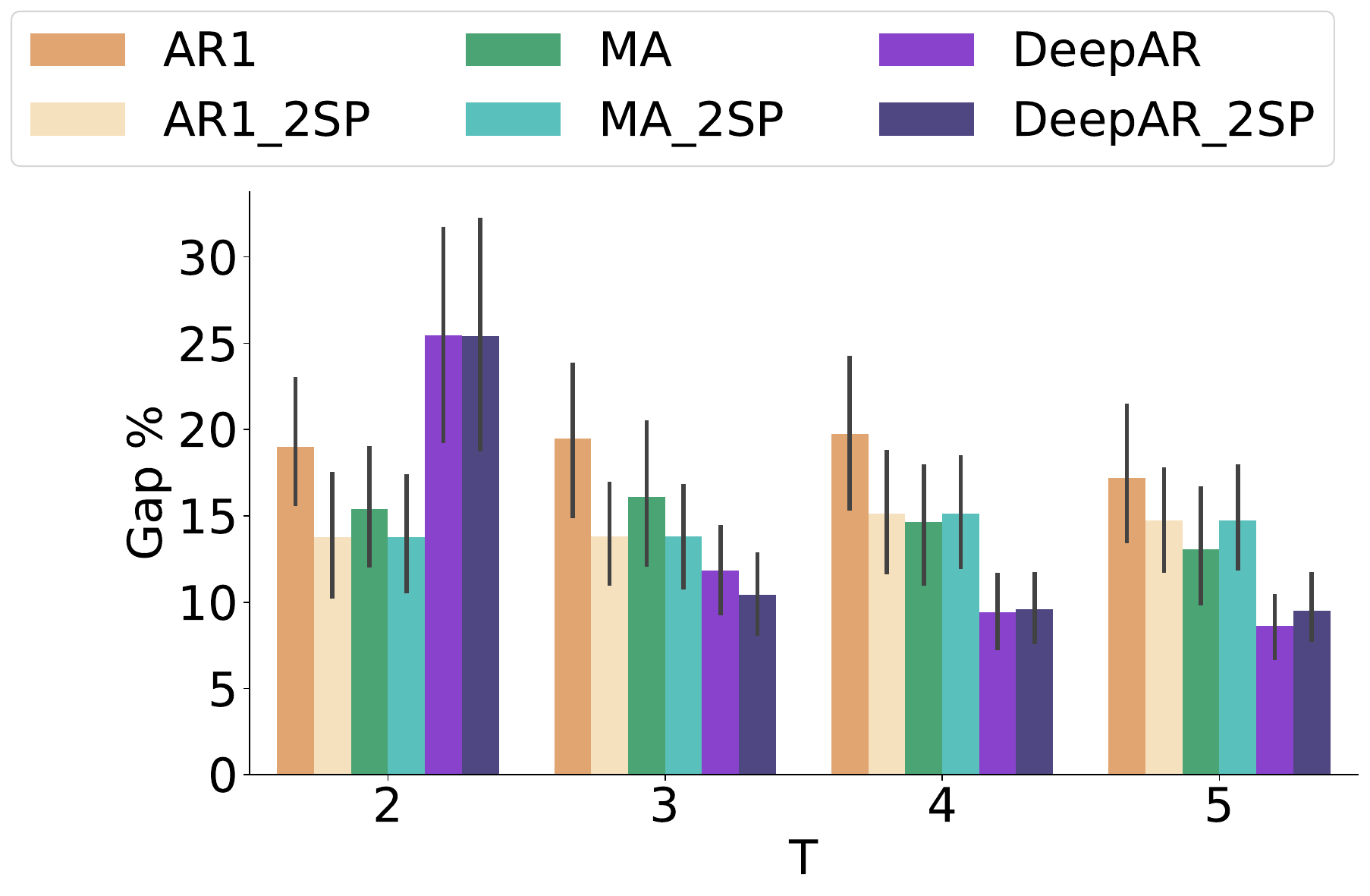}
        \caption{Walmart dataset}
    \end{subfigure}%
    ~ 
    \begin{subfigure}[b]{0.48\textwidth}
        \centering
        \includegraphics[height=1.8in]{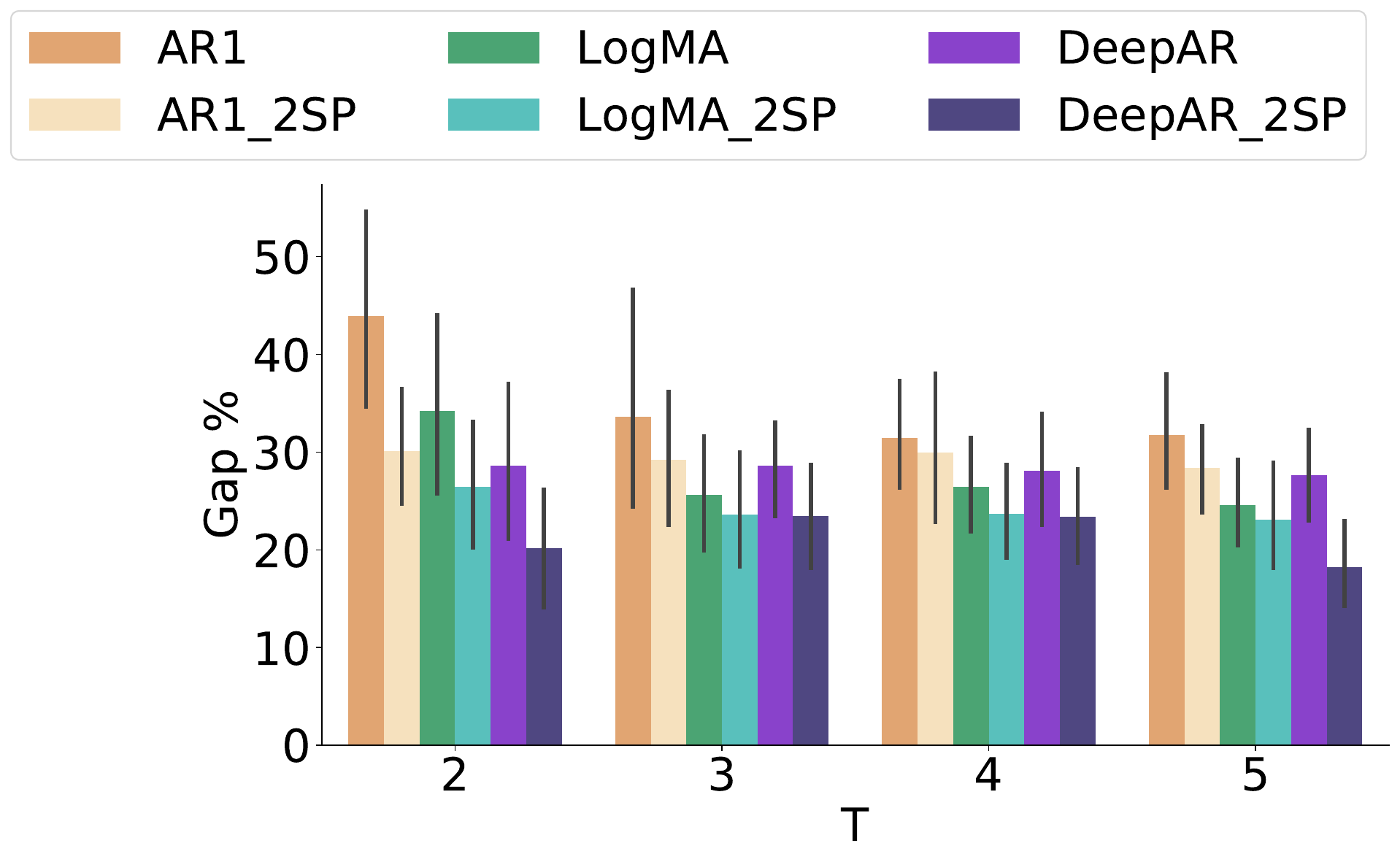}     
        \caption{Favorita dataset}
    \end{subfigure}
    \caption{Gap \% values for deterministic and two-stage policies in the risk-neutral setting}
    \label{fig:stochastic_gapplots}
\end{figure}

Figure~\ref{fig:stochastic_gapplots} provides results on two-stage policies as well as their deterministic counterparts for all the forecasting methods and varying planning horizon lengths. 
First, we observe that, overall, the two-stage policies obtained by using the scenarios sampled from the DeepAR (DeepAR\textunderscore2SP), perform the best. Note that the Walmart dataset with $\Tpar = 2$ is an exception as discussed in our analysis with the deterministic policies and will be omitted in what follows. 

In particular, for the Favorita dataset, two-stage policies consistently provide better results compared to their deterministic counterparts for all the forecasting methods. Most notably, DeepAR\_2SP substantially improves over its deterministic counterpart, and outperforms all the other approaches for this dataset. 

Second, there are various cases where the deterministic policies outperform the stochastic policies. For instance, in the Walmart dataset, DeepAR deterministic policies provide better Gap \% values than ARMA-based two-stage policies. This finding confirms that, well-suited forecasting models for a problem can lead to high-quality solutions via computationally cheap heuristics, while sophisticated optimization models with low-accuracy forecasts are bound to provide low-quality solutions.
Considering the fact that AR(1) model is frequently used in the MSP literature, our study suggests to first assess its predictive performance on the considered application in order not to sacrifice the chance of obtaining high-quality policies.

One important caveat in our assessment is that our evaluations are solely based on a single \red{out-of-sample} scenario, namely, the \red{one consisting of the} true realizations of all the random demand variables. However, in quantifying the value of stochastic solutions, the expected value is used as the measure, that is, the average performance over a (large) set of potential realizations of random variables is taken as the reference, since the ground truth is not known. On the contrary, in our analysis, we intentionally spared some portion of the historical data and evaluated the performance of alternative policies based on the true realizations in order to provide a better picture for decision makers as only one particular scenario is realized in real life. Therefore, although two-stage policies are expected to perform better than the deterministic policies on the average, in our reported results, there are various cases that we observe the opposite as in the case of DeepAR for the Walmart dataset with $\Tpar = \{4,5\}$.

While our analysis with (Log)MA forecasting method for the risk-neutral deterministic policies provide competitive results to those of DeepAR, the gap between the two widens for two-stage policies. As in the deterministic risk-averse case, this can be attributed to the fact that DeepAR provides substantially better probabilistic forecasting results with tight confidence intervals (e.g. see Figure~\ref{fig:forecasting_walmart}, Figure~\ref{fig:forecast_favorita} and $\rho$-risk values provided in Section \ref{sec:forecastingPerformance}), which are used to sample more informative scenarios with a high confidence level.

\red{Lastly, we note that while we considered a risk-averse approach for the deterministic setting using risk-averse forecasts, we have not considered a risk-averse approach for the two-stage setting. A suitable option for the latter can be obtaining a policy based on two-stage robust optimization models, which would be computationally difficult to solve even for small instances in the presence of mixed-integer recourse. 
}
\section{Conclusion}\label{sec:conclusions}
In this paper, we analyze the impact of deep learning-based time-series forecasting on MSP model outcomes. 
In our analysis, we consider two commonly adopted MSP policies (deterministic, and two-stage stochastic), and compare three different forecasting methods, namely, AR(1), (Log)MA, and DeepAR. 
Our numerical study with a practical MSP application provides evidence that improved forecasts can lead to significant benefits in terms of the performance of the look-ahead 
policies. 
Specifically, we observe that, with an improved forecasting method, deterministic heuristics have potential to provide high-quality policies. 
Moreover, once the forecasting quality is assured, the use of stochastic programming can further improve the overall quality of the results. 
We note that these findings align well with those of \citep{gauvin2018stochastic, gauvin2018successive} who emphasize the importance of better forecasting methods for improved policies. 
Lastly, our results suggest that DeepAR can be especially useful as a sampling tool for stochastic programming, as well as to construct uncertainty sets or estimate the worst-case scenario for robust optimization, due to its impressive probabilistic forecasting capabilities.

Considering the success of state-of-the-art \red{deep learning-based} probabilistic 
forecasting methods, as well as their highly suitable features for MSPs, we believe that the proposed approach of combining them with MSP policy generators has significant potential to be utilized widely in practice. As this would allow generating high-quality policies via simpler heuristic, thus substantially reduce the computational burden, it can widen the applicability of multistage stochastic programming. We note that in particular for applications involving complex time series such as stock price prediction, integration of MSP with deep learning-based time series forecasting methods might provide significantly better results than that of with simple forecasting methods that are not able to produce reliable predictions. However, for some other applications such as electricity generation planning where time series shows noticeable patterns such as seasonality and cyclic behaviour which can be accurately predicted with simple autoregressive models, combination of MSP with deep learning-based time-series forecasts might not produce a drastic improvement. 

Our study constitutes the first step in leveraging the power of state-of-the-art probabilistic forecasting methods in stochastic optimization. There are many intriguing future research directions. First, our computational study can be extended through different MSP applications (e.g. hydro-thermal electricity generation) and other probabilistic forecasting methods (e.g. deep state space models~\citep{rangapuram2018deep}), which would be helpful in further assessing the benefits of improved forecasts. Second, in our experiments with two-stage policies, where we used crude Monte Carlo sampling, alternative sampling approaches such as Quasi-Monte Carlo sampling~\citep{homem2014monte} and Monte Carlo Markov Chain based importance sampling~\citep{parpas2015importance} can be employed to see how much variance reduction can be accomplished, and accordingly whether better two-stage policies can be obtained in a reduced computational effort (i.e., using a fewer number of scenarios). Third, it would be interesting to evaluate the impact of improved forecasts on MSP policies constructed via more sophisticated scenario tree-free methods (e.g., decision rule-based approaches) \red{Lastly, the experimental framework that we considered belongs to the ``traditional data-driven sample average approximation approach" \citep{kannan2020data}, as we first constructed the time series forecasting model, and then used the predictions to build the optimization model. In this regard, combining deep learning-based time series forecasting methodologies with other class of modern data-driven optimization frameworks such as ``Predict, then Optimize" \citep{elmachtoub2021smart} is an intriguing future research direction.}

\section*{Acknowledgment}
This research is supported by a grant from LG SciencePark. M. Bodur also acknowledges the support of the Natural Sciences and Engineering Research Council of Canada.



\section*{Data availability statement}
The data that supports the findings of this study are openly available at: \\
\href{https://www.kaggle.com/c/walmart-recruiting-store-sales-forecasting}{\texttt{https://www.kaggle.com/c/walmart-recruiting-store-sales-forecasting}} \\
\href{https://www.kaggle.com/c/favorita-grocery-sales-forecasting}{\texttt{https://www.kaggle.com/c/favorita-grocery-sales-forecasting}}



\bibliographystyle{abbrvnat}
\bibliography{sample}

\newpage

\appendix
\section{ND values for the Walmart dataset}
\label{app:ND_Walmart}
In this section, we provide ND($t',\Tpar$) values for the Walmart dataset when different planning horizon lengths, $T$ values, are assumed. Note that all the ND values are calculated using the same models.
\begin{itemize}
\item $T = 2$
\begin{table}[!ht]
\begin{tabular}{lrr}
\toprule
\textbf{Method}\ /\ $t'$ &     \textbf{0} &     \textbf{1}  \\
\midrule
\textbf{AR(1)}  &  0.22 &  0.23 \\
\textbf{MA}  &  0.20 &  0.22 \\
\textbf{DeepAR} &  0.26 &  0.19 \\
\bottomrule
\end{tabular}
\end{table}

\item $T = 3$ 
\begin{table}[!ht]
\begin{tabular}{lrrr}
\toprule
\textbf{Method}\ /\ $t'$ &     \textbf{0} &     \textbf{1} &     \textbf{2}  \\
\midrule
\textbf{AR(1)}  &  0.27 &  0.31 &  0.26 \\
\textbf{MA}  &  0.24 &  0.27 &  0.22 \\
\textbf{DeepAR} &  0.23 &  0.17 &  0.15 \\
\bottomrule
\end{tabular}
\end{table}

\item $T = 4$ 
\begin{table}[!ht]
\begin{tabular}{lrrrr}
\toprule
\textbf{Method}\ /\ $t'$ &     \textbf{0} &     \textbf{1} &     \textbf{2} &     \textbf{3}  \\
\midrule
\textbf{AR(1)}  &  0.30 &  0.34 &  0.27 &  0.13 \\
\textbf{MA}  &  0.25 &  0.27 &  0.20 &  0.12 \\
\textbf{DeepAR} &  0.20 &  0.15 &  0.13 &  0.11 \\
\bottomrule
\end{tabular}
\end{table} 

\item $T = 5$ 
\begin{table}[!ht]
\begin{tabular}{lrrrrr}
\toprule
\textbf{Method}\ /\ $t'$ &     \textbf{0} &     \textbf{1} &     \textbf{2} &     \textbf{3} &     \textbf{4}  \\
\midrule
\textbf{AR(1)}  &  0.32 &  0.36 &  0.28 &  0.16 &  0.14 \\
\textbf{MA}  &  0.26 &  0.27 &  0.20 &  0.13 &  0.12 \\
\textbf{DeepAR} &  0.19 &  0.15 &  0.13 &  0.12 &  0.13 \\
\bottomrule
\end{tabular}
\end{table} 

\newpage

\item $T = 6$ 
\begin{table}[!ht]
\begin{tabular}{lrrrrrr}
\toprule
\textbf{Method}\ /\ $t'$ &     \textbf{0} &     \textbf{1} &     \textbf{2} &     \textbf{3} &     \textbf{4} &     \textbf{5} \\
\midrule
\textbf{AR(1)}  &  0.32 &  0.35 &  0.26 &  0.16 &  0.12 &  0.07 \\
\textbf{MA}  &  0.25 &  0.25 &  0.19 &  0.13 &  0.12 &  0.09 \\
\textbf{DeepAR} &  0.18 &  0.14 &  0.12 &  0.11 &  0.11 &  0.09 \\
\bottomrule
\end{tabular}
\end{table}

\item $T = 7$ 
\begin{table}[!ht]
\begin{tabular}{lrrrrrrr}
\toprule
\textbf{Method}\ /\ $t'$ &     \textbf{0} &     \textbf{1} &     \textbf{2} &     \textbf{3} &     \textbf{4} &     \textbf{5} &     \textbf{6}  \\
\midrule
\textbf{AR(1)}  &  0.33 &  0.36 &  0.27 &  0.17 &  0.14 &  0.12 &  0.14 \\
\textbf{MA}  &  0.25 &  0.25 &  0.19 &  0.13 &  0.12 &  0.10 &  0.11 \\
\textbf{DeepAR} &  0.17 &  0.14 &  0.12 &  0.12 &  0.12 &  0.11 &  0.13 \\
\bottomrule
\end{tabular}
\end{table}

\end{itemize}

\newpage

\section{ND values for the Favorita dataset}
\label{app:ND_Favorita}
In this section, we provide ND($t',\Tpar$) values for the Favorita dataset when different planning horizon lengths, $T$ values, are assumed. Note that all the ND values are calculated using the same models.
\begin{itemize}

\item $T = 2$
\begin{table}[!ht]
\begin{tabular}{lrr}
\toprule
\textbf{Method}\ /\ $t'$ &     \textbf{0} &     \textbf{1}  \\
\midrule
\textbf{AR(1)}  &  0.86 &  1.58 \\
\textbf{LogMA}  &  0.68 &  0.60 \\
\textbf{DeepAR} &  0.59 &  0.55 \\
\bottomrule
\end{tabular}
\end{table}

\item $T = 3$ 
\begin{table}[!ht]
\begin{tabular}{lrrr}
\toprule
\textbf{Method}\ /\ $t'$ &     \textbf{0} &     \textbf{1} &     \textbf{2} \\
\midrule
\textbf{AR(1)}  &  0.95 &  1.35 &  0.47 \\
\textbf{LogMA}  &  0.61 &  0.51 &  0.42 \\
\textbf{DeepAR} &  0.61 &  0.60 &  0.59 \\
\bottomrule
\end{tabular}
\end{table} 

\item $T = 4$ 
\begin{table}[!ht]
\begin{tabular}{lrrrr}
\toprule
\textbf{Method}\ /\ $t'$ &     \textbf{0} &     \textbf{1} &     \textbf{2} &     \textbf{3}  \\
\midrule
\textbf{AR(1)}  &  0.89 &  1.14 &  0.43 &  0.37 \\
\textbf{LogMA}  &  0.54 &  0.46 &  0.40 &  0.35 \\
\textbf{DeepAR} &  0.53 &  0.49 &  0.44 &  0.30 \\
\bottomrule
\end{tabular}
\end{table}

\item $T = 5$ 
\begin{table}[!ht]
\begin{tabular}{lrrrrr}
\toprule
\textbf{Method}\ /\ $t'$ &     \textbf{0} &     \textbf{1} &     \textbf{2} &     \textbf{3} &     \textbf{4}  \\
\midrule
\textbf{AR(1)}  &  0.89 &  1.11 &  0.44 &  0.40 &  0.42 \\
\textbf{LogMA}  &  0.52 &  0.45 &  0.39 &  0.36 &  0.36 \\
\textbf{DeepAR} &  0.51 &  0.47 &  0.43 &  0.34 &  0.38 \\
\bottomrule
\end{tabular}
\end{table} 

\item $T = 6$ 
\begin{table}[!ht]
\begin{tabular}{lrrrrrr}
\toprule
\textbf{Method}\ /\ $t'$ &     \textbf{0} &     \textbf{1} &     \textbf{2} &     \textbf{3} &     \textbf{4} &     \textbf{5}  \\
\midrule
\textbf{AR(1)}  &  0.93 &  1.16 &  0.49 &  0.46 &  0.53 &  0.61 \\
\textbf{LogMA}  &  0.51 &  0.46 &  0.42 &  0.40 &  0.42 &  0.42 \\
\textbf{DeepAR} &  0.49 &  0.46 &  0.42 &  0.34 &  0.37 &  0.35 \\
\bottomrule
\end{tabular}
\end{table}

\item $T = 7$ 
\begin{table}[!ht]
\begin{tabular}{lrrrrrrr}
\toprule
\textbf{Method}\ /\ $t'$ &     \textbf{0} &     \textbf{1} &     \textbf{2} &     \textbf{3} &     \textbf{4} &     \textbf{5} &     \textbf{6}  \\
\midrule
\textbf{AR(1)}  &  0.97 &  1.20 &  0.52 &  0.51 &  0.58 &  0.63 &  0.47 \\
\textbf{LogMA}  &  0.51 &  0.47 &  0.42 &  0.40 &  0.43 &  0.41 &  0.37 \\
\textbf{DeepAR} &  0.49 &  0.46 &  0.42 &  0.35 &  0.37 &  0.36 &  0.36 \\
\bottomrule
\end{tabular}
\end{table}

\item $T = 8$
\begin{table}[!ht]
\begin{tabular}{lrrrrrrrr}
\toprule
\textbf{Method}\ /\ $t'$ &     \textbf{0} &     \textbf{1} &     \textbf{2} &     \textbf{3} &     \textbf{4} &     \textbf{5} &     \textbf{6} &     \textbf{7}  \\
\midrule
\textbf{AR(1)}  &  1.00 &  1.25 &  0.55 &  0.55 &  0.63 &  0.67 &  0.57 &  0.61 \\
\textbf{LogMA}  &  0.52 &  0.47 &  0.43 &  0.42 &  0.44 &  0.44 &  0.41 &  0.42 \\
\textbf{DeepAR} &  0.48 &  0.45 &  0.41 &  0.35 &  0.38 &  0.37 &  0.36 &  0.36 \\
\bottomrule
\end{tabular}
\end{table}

\item $T = 9$
\begin{table}[!ht]
\begin{tabular}{lrrrrrrrrr}
\toprule
\textbf{Method}\ /\ $t'$ &     \textbf{0} &     \textbf{1} &     \textbf{2} &     \textbf{3} &     \textbf{4} &     \textbf{5} &     \textbf{6} &     \textbf{7} &     \textbf{8}  \\
\midrule
\textbf{AR(1)}  &  1.06 &  1.31 &  0.60 &  0.61 &  0.70 &  0.75 &  0.70 &  0.76 &  0.78 \\
\textbf{LogMA}  &  0.53 &  0.49 &  0.45 &  0.44 &  0.47 &  0.47 &  0.46 &  0.49 &  0.53 \\
\textbf{DeepAR} &  0.48 &  0.46 &  0.42 &  0.36 &  0.39 &  0.39 &  0.39 &  0.41 &  0.46 \\
\bottomrule
\end{tabular}
\end{table}

\newpage

\item $T = 10$
\begin{table}[!ht]
\begin{tabular}{lrrrrrrrrrr}
\toprule
\textbf{Method}\ /\ $t'$ &     \textbf{0} &     \textbf{1} &     \textbf{2} &     \textbf{3} &     \textbf{4} &     \textbf{5} &     \textbf{6} &     \textbf{7} &     \textbf{8} &     \textbf{9} \\
\midrule
\textbf{AR(1)}  &  1.06 &  1.32 &  0.61 &  0.62 &  0.70 &  0.73 &  0.68 &  0.72 &  0.67 &  0.48 \\
\textbf{LogMA}  &  0.53 &  0.50 &  0.46 &  0.45 &  0.47 &  0.48 &  0.47 &  0.49 &  0.49 &  0.42 \\
\textbf{DeepAR} &  0.50 &  0.47 &  0.44 &  0.39 &  0.41 &  0.42 &  0.43 &  0.46 &  0.51 &  0.52 \\
\bottomrule
\end{tabular}
\end{table}

\item $T = 11$
\begin{table}[!ht]
\begin{tabular}{lrrrrrrrrrrr}
\toprule
\textbf{Method}\ /\ $t'$ &     \textbf{0} &     \textbf{1} &     \textbf{2} &     \textbf{3} &     \textbf{4} &     \textbf{5} &     \textbf{6} &     \textbf{7} &     \textbf{8} &     \textbf{9} &     \textbf{10} \\
\midrule
\textbf{AR(1)}  &  1.03 &  1.27 &  0.59 &  0.59 &  0.64 &  0.66 &  0.61 &  0.62 &  0.55 &  0.43 &  0.37 \\
\textbf{LogMA}  &  0.52 &  0.48 &  0.45 &  0.44 &  0.45 &  0.44 &  0.43 &  0.43 &  0.42 &  0.37 &  0.33 \\
\textbf{DeepAR} &  0.48 &  0.46 &  0.43 &  0.38 &  0.40 &  0.40 &  0.40 &  0.41 &  0.43 &  0.40 &  0.28 \\
\bottomrule
\end{tabular}
\end{table}

\end{itemize}

\end{document}